\input ./style/arxiv-general.cfg
\documentclass[aos,MSNbibl,nameyear,dvips]{arximspdf}
\makeatletter
   \@ifpackageloaded{graphicx}{}{\usepackage{graphicx}}
\makeatother
\usepackage{dcolumn}
\usepackage{graphicx}

\doi{10.1214/14-AOS1286}
\volume{43}
\issue{3}
\pubyear{2015}
\firstpage{991}
\lastpage{1026}
\docsubty{FLA}

\makeatletter
\newcolumntype{d}[1]{D{.}{.}{#1}}
\newcommand{\lleft}{\left}
\newcommand{\rrvert}{\vert}
\newcommand{\rright}{\right}
\newcommand{\rrVert}{\Vert}
\newcommand{\llvert}{\vert}
\newcommand{\llVert}{\Vert}
\def\overset{\stackrel}
\newtheorem{theorem}{Theorem}
\newtheorem{corollary}{Corollary}
\newtheorem{proposition}{Proposition}
\newtheorem{lemma}{Lemma}
\newproclaim{example}{Example}
\newproclaim{remark}{Remark}
\newproclaim{definition}{Definition}
\newcommand{\lam}{\lambda}
\newcommand{\argmin}{\mathop{\arg\min}}
\def\Dbar{{\overline{\mathbf{D}}}}
\renewcommand{\epsilon}{\varepsilon}
\makeatother

\begin{document}
\begin{frontmatter}

\title{Asymptotic normality and optimalities in estimation of large
Gaussian graphical models}
\runtitle{Statistical inference for Gaussian graphical model}

\begin{aug}
\author[A]{\fnms{Zhao}~\snm{Ren}\thanksref{m1}\ead[label=e1]{zren@pitt.edu}},
\author[B]{\fnms{Tingni}~\snm{Sun}\thanksref{m2}\ead[label=e2]{tingni@umd.edu}},
\author[C]{\fnms{Cun-Hui}~\snm{Zhang}\thanksref{T2,m3}\ead[label=e3]{cunhui@stat.rutgers.edu}}
\and
\author[D]{\fnms{Harrison~H.}~\snm{Zhou}\corref{}\thanksref{T1,m4}\ead[label=e4]{huibin.zhou@yale.edu}}

\thankstext{T2}{Supported in part by the NSF Grants
DMS-11-06753 and DMS-12-09014 and NSA Grant H98230-11-1-0205.}
\thankstext{T1}{Supported in part by NSF Career Award
DMS-06-45676 and NSF FRG Grant DMS-08-54975.}

\runauthor{Ren, Sun, Zhang and Zhou}

\affiliation{University of Pittsburgh\thanksmark{m1}, University of
Maryland\thanksmark{m2}, Rutgers University\thanksmark{m3} and Yale
University\thanksmark{m4}}

\address[A]{Z. Ren\\
Department of Statistics\\
University of Pittsburgh\\
Pittsburgh, Pennsylvania 15260\\
USA\\
\printead{e1}}

\address[B]{T. Sun\\
Department of Mathematics\\
University of Maryland\\
College Park, Maryland 20742\\
USA\\
\printead{e2}}

\address[C]{C.-H. Zhang\\
Department of Statistics and Biostatistics\\
Hill Center, Busch Campus\\
Rutgers University\\
Piscataway, New Jersey 08854\\
USA\\
\printead{e3}}

\address[D]{H. H. Zhou\\
Department of Statistics\\
Yale University\\
New Haven, Connecticut 06511\\
USA\\
\printead{e4}}
\end{aug}
%

%
\received{\smonth{8} \syear{2013}}
%
\revised{\smonth{10} \syear{2014}}

%
\begin{abstract}
The Gaussian graphical model, a popular paradigm for studying
relationship among variables
in a wide range of applications, has attracted great attention in
recent years.
This paper considers a fundamental question: When is it possible to
estimate low-dimensional parameters at parametric square-root rate in a
large Gaussian graphical model?
A novel regression approach is proposed to obtain asymptotically efficient
estimation of each entry of a precision matrix under a sparseness
condition relative to the sample size.
When the precision matrix is not sufficiently sparse,
or equivalently the sample size is not sufficiently large, a lower
bound is established to show that
it is no longer possible to achieve the parametric rate in the
estimation of each entry.
This lower bound result, which provides an answer to the delicate
sample size question,
is established with a novel construction of a subset of sparse
precision matrices
in an application of Le Cam's lemma. Moreover, the proposed estimator
is proven
to have optimal convergence rate when the parametric rate cannot be
achieved, under
a minimal sample requirement.

The proposed estimator is applied to test the presence of an edge in
the Gaussian graphical
model or to recover the support of the entire model,
to obtain adaptive rate-optimal estimation of the entire precision
matrix as measured
by the matrix $\ell_{q}$ operator norm and to make inference in latent
variables in the graphical model.
All of this is achieved under a sparsity condition on the precision
matrix and a side condition on
the range of its spectrum.
This significantly relaxes the commonly imposed
uniform signal strength condition on the precision matrix,
irrepresentability condition on the Hessian tensor operator of the
covariance matrix
or the $\ell_{1}$ constraint on the precision matrix.
Numerical results confirm our theoretical findings.
The ROC curve of the proposed algorithm, Asymptotic Normal Thresholding (ANT),
for support recovery significantly outperforms that of the popular
GLasso algorithm.
\end{abstract}

%

\end{frontmatter}

\section{Introduction}\label{sec1}

The Gaussian graphical model, a powerful tool for investigating the relationship
among a large number of random variables in a complex system, is used
in a
wide range of scientific applications. A central question for Gaussian
graphical models is how to recover the structure of an undirected
Gaussian graph.
Let $G=(V,E)$ be an undirected graph representing the conditional dependence
relationship between components of a random vector $Z=(Z_{1},\ldots
,Z_{p})^{T}$ as follows. The vertex set $V=\{V_{1},\ldots,V_{p}\}$
represents the components of $Z$. The edge set $E$ consists of pairs $(i,j)$
indicating the conditional dependence between $Z_{i}$ and $Z_{j}$ given all
other components. In applications, the following question is
fundamental: Is
there an edge between $V_{i}$ and $V_{j}$? It is well known that recovering
the structure of an undirected Gaussian graph $G= ( V,E ) $ is
equivalent to recovering the support of the population precision matrix of
the data in the Gaussian graphical model. Let
\[
Z= ( Z_{1},Z_{2},\ldots,Z_{p} ) ^{T}
\sim\mathcal{N} ( \mu,\Sigma ),
\]
where $\Sigma= ( \sigma_{ij} ) $ is the population covariance
matrix. 
The precision matrix, denoted by $\Omega= ( \omega_{ij} )
$, is
defined as the inverse of covariance matrix, $\Omega=\Sigma^{-1}$. There
is an edge between $V_{i}$ and $V_{j}$, that is, $(i,j)\in E$, if and
only if $%
\omega_{ij}\neq0$; see, for example, \citet{GraMod}. Consequently, the
support recovery of the precision matrix $\Omega$ yields the recovery of
the structure of the graph $G$.

Suppose $n$ i.i.d. $p$-variate random vectors $X^{(1)},X^{(2)},\ldots,X^{(n)}$ are observed from the same distribution as $Z$, that is, the
Gaussian $%
\mathcal{N} ( \mu,\Omega^{-1} )$. Assume without loss of
generality that $\mu=0$ hereafter. In this paper, we address the following
two fundamental questions: When is it possible to make statistical inference
for each individual entry of a precision matrix $\Omega$ at the
parametric $%
n^{-1/2}$ rate? When and in what sense is it possible to recover the support
of $\Omega$ in the presence of some small nonzero $|\omega_{ij}|$?

The problems of estimating a large sparse precision matrix and recovering
its support have drawn considerable recent attention. There are mainly two
approaches in the literature. The first approach is a penalized
likelihood estimation
approach with a lasso-type penalty on entries of the precision matrix.
\citet%
{YuanLin07} proposed to use the lasso penalty and studied its asymptotic
properties when $p$ is fixed. \citet{RWRY08} derived the selection consistency
and related error bounds
under an irrepresentability condition on the Hessian tensor
operator and a constraint on the matrix $\ell_{1}$ norm of the precision
matrix. See also \citet{RBLZ08} for Frobenius-based error bounds
and \citet{LFAN09} for concave penalized likelihood estimation without
the irrepresentability condition.
The second approach, proposed earlier by \citet{MB06},
is neighborhood-based. It estimates the precision matrix column by column
by running the lasso or Dantzig selector for each variable against all
the rest of variables; see \citet{Yuan10}, \citet{CLL11}, \citet{CLZ12}
and \citet{SZH122}. The
irrepresentability condition is no longer needed in \citet{CLL11} and
\citet%
{CLZ12} for support recovery, but the thresholding level for support
recovery depends on the matrix $\ell_{1}$ norm of the precision
matrix. The
matrix $\ell_{1}$ norm is unknown and large, which makes the support recovery
procedures there nonadaptive and thus less practical. In \citet{SZH122},
optimal convergence rate in the spectral norm is achieved without requiring
the matrix $\ell_1$ norm constraint or the irrepresentability condition.
However, support recovery properties of the estimator were not analyzed.

In spite of an extensive literature on the topic, the fundamental limit
of support recovery in the Gaussian graphical
model is still largely unknown, let alone an adaptive procedure to
achieve the limit.

Statistical inference of low-dimensional parameters at the $n^{-1/2}$ rate
has been considered in the closely related linear regression model.
\citet%
{SZH12} proposed an efficient scaled lasso estimator of the noise level
under the sample size condition $n\gg(s\log p)^2$, where $s$ is the
$\ell_0$
or capped-$\ell_1$ measure of the size of the unknown regression
coefficient vector.
\citet{zhang2014confidence} proposed an asymptotically normal
low-dimensional projection
estimator (LDPE) for the regression coefficients under the same sample
size condition.
Both estimators converge at the $n^{-1/2}$ rate, and their
asymptotic efficiency can be understood from the
minimum Fisher information in a more general context [\citet{ZhangOber11}].
A proof of the asymptotic efficiency of the LDPE
was given in \citet{GeerBR13} where the generalized linear
model was also considered.
Alternative methods for testing and estimation of regression
coefficients were proposed
in \citet{BelloniCH12}, \citet{Buhlmann12}, \citet{JavanmardM13} and
\citet{Liu13}.
However, the optimal rate of convergence is unclear from these papers
when the sample size
condition $n\gg(s\log p)^2$ fails to hold.
Please see Section~\ref{Discu1} for more details of their connection
with this paper.

This paper makes important advancements in the understanding of statistical
inference of low-dimensional parameters in the Gaussian graphical model in
the following ways. Let $s$ be the maximum degree of the graph or a certain
more relaxed capped-$\ell_1$ measure of the complexity of the precision
matrix. We prove that the estimation of each $\omega_{ij}$ at the parametric
$n^{-1/2}$ convergence rate requires the sparsity condition $s=
O(n^{1/2}/\log p)$ or equivalently a sample size of order $(s\log
p)^2$. We
propose an adaptive estimator of individual $\omega_{ij}$ and prove its
asymptotic normality and efficiency when $n \gg(s\log p)^2$. Moreover, we
prove that the proposed estimator achieves the optimal convergence rate when
the sparsity condition is relaxed to $s\le c_0 n/\log p$ for a certain
positive constant $c_0$. The efficient estimator of the individual $%
\omega_{ij}$ is then used to construct fully data-driven procedures to
recover the support of $\Omega$ and to make statistical inference about
latent variables in the graphical model.

The methodology we are proposing is a novel regression approach briefly
described in \citet{SunZhang2012SinicaComment}. In this regression approach,
the main task is not to estimate the slope, as seen in \citet{MB06},
\citet%
{Yuan10}, \citet{CLL11}, \citet{CLZ12} and \citet{SZH12}, but to
estimate the
noise level. For a vector $Z$ of length $p$ and any index subset $A$ of
$ \{ 1,2,\ldots,p \}$, we denote by $Z_{A}$ the sub-vector
of $Z$
with elements indexed by $A$. Similarly for a
matrix $U$ and two index subsets $A$ and $B$ of
$ \{ 1,2,\ldots,p \}$, we denote by $U_{A,B}$ the
$|A|\times|B|$
sub-matrix of $U$ with elements in rows in $A$ and columns in $B$.
Consider $A= \{ i,j \}$ with $i\neq j$, so that $Z_{A}= (
Z_{i},Z_{j} ) ^{T}$ and $\Omega_{A,A}=\bigl(
{{\omega_{ii}\atop  \omega_{ji}} \enskip{ \omega_{ij} \atop \omega_{jj}}}
 \bigr)$. It is well known that
\[
Z_{A}|Z_{A^{c}}\sim\mathcal{N} \bigl( -\Omega_{A,A}^{-1}
\Omega _{A,A^{c}}Z_{A^{c}},\Omega_{A,A}^{-1} \bigr).
\]
This observation motivates us to consider the estimation of individual
entries of $\Omega$,
$\omega_{ii}$ and $\omega_{ij}$,
by estimating the noise level in the regression of the two response
variables in $A$
against the variables in $A^c$.
The noise level $\Omega_{A,A}^{-1}$ has only three parameters.
When $\Omega$ is sufficiently sparse, a penalized regression
approach is proposed in Section~\ref{methodinf} to obtain an asymptotically
efficient estimation of $\omega_{ij}$ in the following sense: The
estimator is
asymptotically normal, and its asymptotic variance matches that of the maximum
likelihood estimator in the classical setting where the dimension $p$
is a
fixed constant. Consider the class of parameter spaces modeling sparse
precision matrices with at most $k_{n,p}$ nonzero elements in
each column,
%
\begin{equation}
\label{sparseparaspace} \mathcal{G}_{0}(M,k_{n,p})=\lleft\{
\begin{array} {c}\displaystyle  \Omega= ( \omega_{ij} ) _{1\leq i,j\leq p}
\dvtx \max_{1\leq
j\leq
p}\sum_{i=1}^p
1 \{ \omega_{ij}\neq0 \} \leq k_{n,p},
\\
\mbox{and }\displaystyle 1/M\leq\lambda_{\min} ( \Omega ) \leq \lambda _{\max
}
( \Omega ) \leq M%
\end{array} %
 \rright\} ,
\end{equation}
where $1 \{ \cdot \} $ is the indicator function, and $M$
is some
constant greater than $1$. The following theorem shows that a necessary and
sufficient condition to obtain a $n^{-1/2}$-consistent estimation of
$\omega
_{ij}$ is $k_{n,p}=O ( \frac{\sqrt{n}}{\log p} ) $, and
when $%
k_{n,p}=o ( \frac{\sqrt{n}}{\log p} ) $ the procedure to be
proposed in Section~\ref{methodinf} is asymptotically efficient.

\begin{theorem}
\label{Mainl0} Let $X^{(i)}{\overset{\mathit{i.i.d.}}{\sim}}\mathcal
{N}%
_{p}(\mu,\Sigma)$, $i=1,2,\ldots,n$. Assume that $3 \leq k_{n,p}\leq
c_{0}n/\log
p$ with a sufficiently small constant $c_{0}>0$ and $p\geq k_{n,p}^{\nu}$
with some $\nu>2$. 
\begin{longlist}[(ii)]
\item[(i)] There exists a constant $\epsilon_{0}>0$ such that
\[
\inf_{i,j} \inf_{\hat{\omega}_{ij}}\sup
_{\mathcal
{G}_{0}(M,k_{n,p})}\mathbb{%
P} \bigl\{\llvert \hat{
\omega}_{ij}-\omega_{ij}\rrvert \geq \epsilon _{0}
\max \bigl\{n^{-1}k_{n,p}\log p,n^{-1/2} \bigr\} \bigr\}
\geq \epsilon_{0}.
\]
Moreover, the minimax risk of estimating $\omega_{ij}$ over the class $
\mathcal{G}_{0}(M,k_{n,p})$ satisfies
%
\begin{equation}
\inf_{\hat{\omega}_{ij}}\sup_{\mathcal{G}_{0}(M,k_{n,p})}\mathbb {E}%
\llvert \hat{\omega}_{ij}-\omega_{ij}\rrvert \asymp\max \bigl
\{n^{-1}k_{n,p}\log p,n^{-1/2} \bigr\} \label{Minimaxelement}
\end{equation}
uniformly in $(i,j)$, provided that $n=O ( p^{\xi} ) $
with some $\xi>0$.
\item[(ii)] The estimator $\hat{\omega}_{ij}$ defined in (\ref
{estimator}) in Section~\ref{Methodology} 
is rate optimal in the sense of
\[
\lim_{(C,n)\rightarrow(\infty,\infty)}\max_{i,j}\sup
_{\mathcal
{G}_{0}(M,k_{n,p})}\mathbb{P} \bigl\{\llvert \hat{%
\omega}_{ij}-\omega_{ij}\rrvert \geq C\max \bigl\{
n^{-1}k_{n,p}\log p,n^{-1/2} \bigr\} \bigr\} = 0.
\]
%
Furthermore, the estimator $\hat{\omega}_{ij}$ 
is asymptotically efficient when $k_{n,p}=o ( \frac{\sqrt
{n}}{\log
p}%
 )$, that is, with $F_{ij}=(\omega_{ii}\omega_{jj}+\omega
_{ij}^{2})^{-1}$
being the Fisher information for estimating $\omega_{ij}$ and
${\hat F}_{ij}=({\hat\omega}_{ii}{\hat\omega}_{jj}+{\hat\omega
}_{ij}^{2})^{-1}$ its estimate,
%
\begin{equation}
\sqrt{n{\hat F}_{ij}} ( \hat{\omega}_{ij}-
\omega_{ij} ) \overset{D} {%
\rightarrow}\mathcal{N} ( 0,1 ),\qquad {
\hat F}_{ij}/F_{ij}\to 1. \label{result 4}
\end{equation}
\end{longlist}
\end{theorem}

The lower bound is established through Le Cam's lemma and a novel
construction of a subset of sparse precision matrices. An important
implication of the lower bound is that the difficulty of support recovery
for sparse precision matrices is different from that for sparse covariance
matrices when\vspace*{-1pt} $k_{n,p}\gg ( \frac{\sqrt{n}}{\log p} ) $, and
when $%
k_{n,p}=o ( \frac{\sqrt{n}}{\log p} ) $ the difficulty of support
recovery for sparse precision matrices is just the same as that for sparse
covariance matrices.

It is worthwhile to point out that the asymptotic efficiency result is
obtained without the need to assume the irrepresentability condition or
the $%
\ell_{1}$ constraint of the precision matrix which are commonly
required in
the literature. For preconceived $(i,j)$, two immediate consequences of
(\ref{result 4}) are efficient interval estimation of $\omega_{ij}$ and
efficient test for the existence of an edge between $V_{i}$
and $V_{j}$ in the graphical model, that is, the hypotheses $\omega_{ij}=0$.
However, the impact of Theorem~\ref{Mainl0} is much broader.
We derive fully adaptive thresholded versions of the estimator and
prove that the thresholded estimators achieve rate optimality in support
recovery without assuming the irrepresentability condition and
in various matrix norms for the estimation of the entire precision
matrix $\Omega$
under weaker assumptions than the requirements of existing results in
the literature.
In addition, we extend our inference and estimation framework to a
class of latent
variable graphical models. See Section~\ref{app} for details.

Our work on optimal estimation of precision matrices given in the present
paper is closely connected to a growing literature on the estimation of large
covariance matrices. Many regularization methods have been proposed and
studied. For example, \citeauthor{BL08A} (\citeyear{BL08A,BL08B}) proposed banding and thresholding
estimators for estimating bandable and sparse covariance matrices,
respectively, and obtained rate of convergence for the two estimators. See
also \citet{ELK08} and \citet{LFAN09}. \citet{CZZH10} established the optimal
rates of convergence for estimating bandable covariance matrices. \citet
{CZh12Sparse} and \citet{CLZ12} obtained the minimax rate of
convergence for
estimating sparse covariance and precision matrices under a range of losses
including the spectral norm loss. In particular, a new general lower bound
technique for matrix estimation was developed there.
More recently, \citet{SZH122} proposed to apply a scaled lasso to estimate
$\Omega$ and proved its rate optimality in spectrum norm without imposing
an $\ell_1$ norm assumption on $\Omega$.

The proposed estimator was briefly described in \citet%
{SunZhang2012SinicaComment} along with a statement of the efficiency of the
estimator without proof under the sparsity assumption
$k_{n,p}=o(n^{-1/2}\log p)$. While we are working on the delicate issue
of the necessity
of the sparsity condition $k_{n,p}=o(n^{1/2}/\log p)$ and the
optimality of
the method for support recovery and estimation under the general sparsity
condition $k_{n,p}=o(n/\log p)$, \citet{Liu13} developed $p$-values for
testing $\omega_{ij}=0$ and related FDR control methods under the stronger
sparsity condition $k_{n,p}=o(n^{1/2}/\log p)$. However, his method cannot
be directly converted into confidence intervals, and the optimality of
his method is
unclear under either sparsity conditions.

The paper is organized as follows. In Section~\ref{methodinf}, we introduce
our methodology and main results for statistical inference.
Applications to
the estimation under the spectral norm, support recovery and the
estimation of
latent variable graphical models are presented in Section~\ref{app}.
Results on linear regression are presented in Section~\ref{regression-revisit}
to support the main theory.
Section~\ref{Discu0} discusses possible extensions of our results and the
connection between our and existing results. Numerical studies are
presented in Section~\ref{Numer}.
The proof for the novel lower bound result is given in Section~\ref{Thmproof5}. Additional proofs 
are provided in \citet{Ren13supp}.

\textit{Notation}. We summarize here some notation to be used throughout
the paper. For $1\leq w\leq\infty$, we use $\llVert  u\rrVert _{w}$
and $\llVert  A\rrVert _{w}$ to denote the usual vector $\ell
_{w}$ norm,
given a vector $u\in\mathbb{R}^{p}$ and a matrix $A= (
a_{ij} )
_{p\times p}$, respectively. In particular, $\llVert  A\rrVert
_{\infty
}$ denote the entry-wise maximum $\max_{ij}\llvert  a_{ij}\rrvert
$.
We shall write $\llVert \cdot\rrVert $ without a subscript
for the
vector $\ell_{2}$ norm. The matrix $\ell_{w}$ operator norm of a
matrix $A$ is
defined by $|\!|\!|A|\!|\!|_{w}=\max_{\Vert x\Vert_{w}=1}\Vert Ax\Vert_{w}$. The
commonly used spectral norm $|\!|\!|\cdot|\!|\!|$ coincides with the matrix
$\ell
_{2}$ operator norm $|\!|\!|\cdot|\!|\!|_{2}$.

\section{Methodology and statistical inference}
\label{Methodology}
\label{methodinf}

In this section we introduce our methodology for estimating each entry
and more generally, a smooth functional of any square submatrix of fixed
size. Asymptotic efficiency results are stated in Section~\ref{infer.sec}
under a sparseness assumption. The lower bound in Section~\ref{low.sec}
shows that the sparseness condition is sharp for the asymptotic efficiency
proved in Section~\ref{infer.sec}.

\subsection{Methodology}

We will first introduce the methodology to estimate each entry $\omega
_{ij}$%
, and discuss its extension to the estimation of functionals of a submatrix
of the precision matrix.

The methodology is motivated by the following simple observation with $%
A= \{ i,j \} $:
%
\begin{equation}
Z_{ \{ i,j \} }|Z_{ \{ i,j \} ^{c}}\sim\mathcal {N} \bigl( -
\Omega_{A,A}^{-1}\Omega_{A,A^{c}}Z_{ \{ i,j \}
^{c}},\Omega
_{A,A}^{-1} \bigr). \label{CondiDis2}
\end{equation}
Equivalently we write a bivariate linear model
%
\begin{equation}
( Z_{i},Z_{j} ) =Z_{ \{ i,j \} ^{c}}^{T}{
\beta}%
+ ( \eta_{i},\eta_{j} ) , \label{regression}
\end{equation}
where the coefficients and error distributions are
%
\begin{equation}
\bolds{\beta}=\bolds{\beta}_{A^c,A} = -\Omega_{A^{c},A}\Omega
_{A,A}^{-1}, \qquad ( \eta _{i},\eta_{j} )
^{T}\sim\mathcal{N} \bigl( 0,\Omega _{A,A}^{-1}
\bigr) . \label{true beta}
\end{equation}
Denote the covariance matrix of $ ( \eta_{i},\eta_{j} )
^{T}$ by
\[
\Theta_{A,A}=\Omega_{A,A}^{-1}=\pmatrix{ \theta_{ii} & \theta_{ij}
\vspace*{2pt}\cr
\theta_{ji} & \theta_{jj}}.
\]
We will estimate $\Theta_{A,A}$ and expect that an efficient estimator
of $%
\Theta_{A,A}$ yields an efficient estimation of the entries of $\Omega_{A,A}$
by inverting the estimator of $\Theta_{A,A}$.

Denote the $n$ by $p$-dimensional data matrix by $\mathbf{X.}$ The
$i$th row
of the data matrix is the $i$th sample $X^{(i)}$. Let $\mathbf{X}_{A}$
be the
sub-matrix of $\mathbf{X}$ composed of columns indexed by $A$. Based on
the regression
interpretation (\ref{regression}), we have the following data version
of the
multivariate regression model
%
\begin{equation}
\mathbf{X}_{A}=\mathbf{X}_{A^{c}}\bolds{\beta}+\bolds{
\epsilon}_{A}. \label{regression-data}
\end{equation}
Here each row of (\ref{regression-data}) is a sample of the linear
model (%
\ref{regression}). Note that $\bolds{\beta}=\bolds{\beta}_{A^c,A}$
is a $p-2$ by $2$-dimensional
coefficient matrix. Denote a sample version of $\Theta_{A,A}$ by
%
\begin{equation}
\Theta_{A,A}^{\mathrm{ora}}= \bigl( \theta_{ij}^{\mathrm{ora}}
\bigr) _{i\in A,j\in A}= 
\bolds{\epsilon}_{A}^{T}
\bolds{\epsilon}_{A}/n, \label{oracle}
\end{equation}
which is an oracle MLE of $\Theta_{A,A}$ based on the extra knowledge
of $\bolds{%
\beta}$. The oracle MLE of $\Omega_{A,A}$ is
%
\begin{equation}
\Omega_{A,A}^{\mathrm{ora}}= \bigl( \omega_{ij}^{\mathrm{ora}}
\bigr) _{i\in A,j\in
A}= \bigl( \Theta_{A,A}^{\mathrm{ora}} \bigr)
^{-1}. \label{oracle1}
\end{equation}
Of course $\bolds{\beta}$ is unknown, and we will need to estimate $%
\bolds{\beta}$ and plug in its estimator to estimate $\bolds{\epsilon}%
_{A}$. This general scheme can be formally written as
%
\begin{equation}
\hat{\Omega}_{A,A}= ( {\hat\omega}_{ij} ) _{i\in A,j\in
A}=
\hat {\Theta}_{A,A}^{-1},\qquad 
\hat{\Theta}_{A,A}=
( \hat{\theta}_{ij} ) _{i,j\in A} = \hat{\bolds{
\epsilon}}_{A}^{T}\hat{\bolds {\epsilon}}%
_{A}/n,
\label{estimator}
\end{equation}
where $\hat{\bolds {\epsilon}}_{A}$ is the estimated residual
corresponding to a suitable
estimator of $\bolds{\beta}_{A^c,A}$, that is,
%
\begin{equation}
\hat{\bolds {\epsilon}}_{A}=\mathbf{X}_{A}-\mathbf
{X}_{A^{c}}\hat{\bolds{\beta}}_{A^c,A}. \label{estimator1}
\end{equation}

Now we introduce specific estimators of $\hat{\bolds {\beta
}}=\hat{\bolds{\beta}}_{A^c,A}=(\hat{\bolds {\beta}}_{i},\hat{\bolds{\beta}}_{j})$.
For each $m\in A= \{i,j \} $, we apply a scaled lasso
estimator to the univariate linear
regression of $\mathbf{X}_{m}$ against $\mathbf{X}_{A^{c}}$ as follows:
%
\begin{equation}
\bigl\{ \hat{\bolds {\beta}}_{m},\hat{\theta}_{mm}^{1/2}
\bigr\} =\argmin_{b\in\mathbb{R}^{p-2},\sigma\in\mathbb{R}^{+}} \biggl\{ \frac{
\llVert \mathbf{X}_{m}-\mathbf{X}_{A^{c}}b\rrVert
^{2}}{2n\sigma
}+%
\frac{\sigma}{2}+\lambda\sum_{k\in A^{c}}\frac{\llVert \mathbf
{X}%
_{k}\rrVert }{\sqrt{n}}
\llvert b_{k}\rrvert \biggr\} , \label{Scaled Lasso}
\end{equation}
with a weighted $\ell_{1}$ penalty, where the vector $b$ is indexed by $A^{c}$.
This is equivalent to standardizing the design vector to length $\sqrt
{n}$ and then
applying the $\ell_1$ penalty to the new coefficients $(\llVert
\mathbf{X}_{k}\rrVert /\sqrt{n})b_{k}$.
The penalty level $\lam$ will be specified explicitly later.
It can be shown that the definition of $\hat{\theta}_{mm}$ in (\ref
{estimator}) is consistent
with the $\hat{\theta}_{mm}$ obtained from the scaled lasso (\ref
{Scaled Lasso})
for each $m\in A$ and each $A$. We also consider the following least
squares estimator (LSE)
in the model ${\hat S}_{mm}$ selected in (\ref{Scaled Lasso}):
%
\begin{equation}
\bigl\{ \hat{\bolds {\beta}}_{m},\hat{\theta}_{mm}^{1/2}
\bigr\} =\argmin_{b\in\mathbb{R}^{p-2},\sigma\in\mathbb{R}^{+}} \biggl\{ \frac{
\llVert \mathbf{X}_{m}-\mathbf{X}_{A^{c}}b\rrVert
^{2}}{2n\sigma
}+%
\frac{\sigma}{2}\dvtx \operatorname{ supp}(b) \subseteq{\hat S}_{mm}
\biggr\} , \label{Scaled Lasso LSE}
\end{equation}
where $\operatorname{ supp}(b)$ denotes the support of vector $b$.

Different versions of scaled lasso, in the sense of scale-free
simultaneous estimation
of the regression coefficients and noise level, have been considered in
\citet{StadlerBG10}, \citet{Antoniadis10} and
\citeauthor{SunZ10} (\citeyear{SunZ10,SZH12})
among others.
The $\hat{\bolds {\beta}}_{m}$ in (\ref{Scaled Lasso})
is equivalent to the square-root lasso in \citet{BelloniCW11}.
Theoretical properties of the LSE after model selection, given in
(\ref{Scaled Lasso LSE}), were studied in \citeauthor{SZH12}
(\citeyear{SZH12,SZH122}).


Our methodology can be routinely extended into a more general form. For any
subset $B\subset \{ 1,2,\ldots,p \} $ with a bounded size, the
conditional distribution of $Z_{B}$ given $Z_{B^{c}}$ is
%
\begin{equation}
Z_{B}|Z_{B^{c}}=\mathcal{N} \bigl( -\Omega_{B,B}^{-1}
\Omega _{B,B^{c}}Z_{B^{c}},\Omega_{B,B}^{-1}
\bigr) , \label{General conditional dist}
\end{equation}
so that the associated multivariate linear regression model is $\mathbf
{X}%
_{B}=\break \mathbf{X}_{B^{c}}\bolds{\beta}_{B,B^{c}}+\bolds{\epsilon}_{B}$
with $\bolds{\beta}_{B^{c},B}=-\Omega_{B^{c},B}\Omega_{B,B}^{-1}$
and $%
\epsilon_B\sim\mathcal{N} ( 0,\Omega_{B,B}^{-1} )$.
Consider a
more general problem of estimating a smooth functional of $\Omega
_{B,B}^{-1}$, denoted by
\[
\zeta=\zeta \bigl( \Omega_{B,B}^{-1} \bigr).
\]
When $\bolds{\beta}_{B^{c},B}$ is known, $\epsilon_B$ is sufficient
for $%
\Omega_{B,B}^{-1}$ due to the independence of $\epsilon_B$ and
$\mathbf
{X}%
_{B^{c}}$, so that an oracle maximum likelihood estimator of $\zeta$
can be defined as
\[
\zeta^{\mathrm{ora}}=\zeta \bigl( \bolds{\epsilon}_{B}^{T}
\bolds{\epsilon }%
_{B}/n \bigr).
\]
We apply an adaptive regularized estimator $\hat{\bolds {\beta}}_{B^{c},B}$
by regressing $\mathbf{X}_B$ against $\mathbf{X}_{B^{c}}$, for example,
a penalized LSE or the LSE after model selection.
We estimate the residual matrix by $\hat{\bolds {\epsilon
}}_{B}=\mathbf
{X}_{B}-\mathbf{X}%
_{B^{c}}\hat{\bolds {\beta}}_{B^{c},B}$, and $%
\zeta ( \Omega_{B,B}^{-1} ) $ by
%
\begin{equation}
\hat{\zeta}=\zeta \bigl( \hat{\bolds {\epsilon}}_{B}^{T}
\hat{\bolds{\epsilon%
}}_{B}/n \bigr). \label{estimator3}
\end{equation}

\subsection{Computational complexity}
For statistical inference about a single entry $\omega_{ij}$ of the
precision matrix $\Omega$
with preconceived $i$ and $j$, the computational cost of the estimator
(\ref{estimator}) is of
the same order as that of a single run of the scaled lasso~(\ref
{Scaled Lasso}).

For the estimation of the entire precision matrix $\Omega$, the
definition of (\ref{estimator}) requires
the computation of ${\hat\omega}_{ij}$ for ${p\choose2}$ different
$A=\{i,j\}$, $i<j$.
However, the computational cost for these ${p\choose2}$ different
${\hat\omega}_{ij}$
is no greater than that of $(1+{\bar s})p$ runs of (\ref{Scaled
Lasso}) where
${\bar s}$ is the average size of the selected model for regressing a
single $\mathbf{X}_j$ against
the other $p-1$ variables. This can be seen as follows.
Define the ``one-versus-rest'' estimator as
\[
\Bigl\{ \hat{\bolds {\beta}}^{(1)}_{-i,i},\sqrt{\hat{
\theta }_{ii}^{(1)}} \Bigr\} = \argmin_{b\in\mathbb{R}^{p-1},\sigma\in\mathbb{R}%
^{+}} \biggl\{
\frac{\llVert \mathbf{X}_{i}-\mathbf{X}_{\{i\}
^c}b\rrVert
^{2}}{2n\sigma}+\frac{\sigma}{2}+\lambda\sum_{k\neq i}
\frac{\llVert
\mathbf{X}_{k}\rrVert }{\sqrt{n}}\llvert b_{k}\rrvert \biggr\}
\]
and ${\hat S}_i^{(1)}= \operatorname{supp}(\hat{\bolds {\beta}}_{-i,i}^{(1)})$.
For $j\notin\{i\}\cup{\hat S}_i^{(1)}$, the ``two-versus-rest''
estimator (\ref{Scaled Lasso}) satisfies
$ \{ \hat{\bolds {\beta}}_{m},\hat{\theta}_{mm}^{1/2} \}
=  \{ \hat{\bolds {\beta}}^{(1)}_{\{i,j\}^c,i},\sqrt{\hat
{\theta
}_{ii}^{(1)}} \}$
when $m=i$ and $A=\{i,j\}$. Thus we only need to carry out $1+|{\hat
S}_i^{(1)}|$
runs of (\ref{Scaled Lasso}) to compute the two-versus-rest estimator
$ \{ \hat{\bolds {\beta}}_{m},\hat{\theta}_{mm}^{1/2} \}$
for all $m=i$ and $A=\{i,j\}$, $j\neq i$, where $|{\hat S}_i^{(1)}|$
denotes the cardinality of the set ${\hat S}_i^{(1)}$. Consequently,
the total required runs of
the scaled lasso (\ref{Scaled Lasso}) is
$\sum_{i=1}^p (1+|{\hat S}_i^{(1)}|) = (1+{\bar s})p$.
It follows from Theorem~\ref{Error bound of LSE} below
that $(1+{\bar s})p$ is of the order $\#\{(i,j)\dvtx \omega_{ij}\neq0\}$.
Thus for the computation of the
estimator (\ref{estimator}) for the entire precision matrix $\Omega$,
the order of the total number of runs of (\ref{Scaled Lasso}) is
the total number of edges of the graphical model corresponding to
$\Omega$.

\subsection{Statistical inference}
\label{infer.sec}
Our analysis can be outlined as follows.
We prove that estimators in the form of (\ref{estimator})
possess the asymptotic normality and efficiency properties
claimed in Theorem~\ref{Mainl0} when the following conditions
hold for certain fixed constant $C_0$, $\epsilon_{\Omega}\to0$ and all
$\delta\ge1$:
%
\begin{eqnarray}
\label{cond-1} \max_{A:A=\{i,j\}} \mathbb{P} \bigl\{\bigl\llVert
\mathbf{X}_{A^{c}}(\hat{\bolds {\beta }}_{A^c,A} -
\beta_{A^c,A})\bigr\rrVert ^2 \ge C_0 s \delta
\log p \bigr\}&\le& p^{-\delta+1}\epsilon_\Omega,
\\
\label{cond-2}  \max_{A:A=\{i,j\}} \mathbb{P} \bigl\{\bigl\llVert
\Dbar _{A^c}^{1/2}(\hat{\bolds{\beta}}_{A^c,A} -
\beta_{A^c,A})\bigr\rrVert _1 \ge C_0s\sqrt{
\delta(\log p)/n} \bigr\}& \le& p^{-\delta+1}\epsilon _\Omega,
\end{eqnarray}
with $\Dbar= \operatorname{ diag}(\mathbf{X}^T\mathbf{X}/n)$,
and for $\theta_{ii}^{\mathrm{ora}}=\|\mathbf{X}_{i} - \mathbf
{X}_{A^{c}}\bolds{{\beta}}_{A^c,i}\|^2/n$,
%
\begin{equation}
\label{cond-3} \max_{A:A=\{i,j\}} \mathbb{P} \biggl\{\biggl |
\frac{{\hat\theta}_{ii}}{
\theta_{ii}^{\mathrm{ora}}}-1 \biggr| \ge C_0 s\delta(\log p)/n \biggr\}\le
p^{-\delta+1}\epsilon_\Omega, 
\end{equation}
with a certain complexity measure $s$ of the precision matrix $\Omega$,
provided that the spectrum of $\Omega$ is bounded, and the sample size
$n$ is no smaller than
$(s\log p)^2/c_0$ for a sufficiently small $c_0>0$.
This is carried out by comparing the estimator in~(\ref{estimator})
with the oracle
MLE in (\ref{oracle}) and (\ref{oracle1}) and proving
\[
\kappa_{ij}^{\mathrm{ora}}=\sqrt{n}\frac{\omega_{ij}^{\mathrm{ora}}-\omega
_{ij}}{%
\sqrt{\omega_{ii}\omega_{jj}+\omega_{ij}^{2}}} \to\mathcal {N}
( 0,1 ),
\]
or equivalently the asymptotic normality of the oracle
MLE in (\ref{oracle1}) with mean $\omega_{ij}$ and variance
$n^{-1} ( \omega_{ii}\omega_{jj}+\omega_{ij}^{2} ) $.
We then prove (\ref{cond-1}), (\ref{cond-2}) and (\ref{cond-3}) for both
the scaled lasso estimator (\ref{Scaled Lasso}) and the LSE after the
scaled lasso selection
(\ref{Scaled Lasso LSE}).
Moreover, we prove that certain thresholded versions of the proposed estimator
possesses global optimality properties, as discussed below Theorem~\ref
{Mainl0}, under
the same boundedness condition on the spectrum of $\Omega$ and
a more relaxed condition on the sample size.

For the $\ell_0$ class $\mathcal{G}_{0}(M,k_{n,p})$ in (\ref
{sparseparaspace}),
the complexity measure for the precision matrix $\Omega$ is the maximum
degree $s=k_{n,p}$
of the corresponding graph.
The $\ell_0$ complexity measure can be relaxed to a capped-$\ell_1$
measure as follows.
For $\lambda>0$, define capped-$\ell_1$ balls as
%
\begin{equation}
\mathcal{G}^*(M,s,\lambda)= \bigl\{ \Omega\dvtx s_\lam(\Omega)\le s,
1/M\leq \lambda_{\min} ( \Omega ) \leq\lambda_{\max} ( \Omega )
\leq M \bigr\}, \label{StarBall}
\end{equation}
where $s_\lam=s_{\lambda}(\Omega)=\max_{j}\sum_{i=1}^p\min \{
1,\llvert \omega_{ij}\rrvert /\lambda \}$ for $\Omega
= (\omega_{ij} ) _{1\leq i,j\leq p}$. In this paper,
$\lambda
$ is
of the order $\sqrt{(\log p)/n}$. We omit the subscript $\lambda$
from $s$
when it is clear from the context. When $\llvert \omega_{ij}\rrvert $
is either $0$ or larger than $\lambda$, $s_{\lambda}$ is the maximum node
degree of the graph. In general, maximum node degree is an upper bound of
the capped-$\ell_1$ measure $s_{\lambda}$. The spectrum of $\Sigma$ is
bounded in the matrix class $\mathcal{G}^{\ast}(M,s,\lambda)$ as in
the $%
\ell_{0}$ ball (\ref{sparseparaspace}). 
The following theorem bounds the difference between our estimator and
the oracle
estimators and the difference between the standardized oracle estimator
and a standard normal variable.

\begin{theorem}
\label{Main} Let $\Theta_{A,A}^{\mathrm{ora}}$ and $\Omega_{A,A}^{\mathrm{ora}}$ be the
oracle MLE
defined in (\ref{oracle}) and (\ref{oracle1}), respectively, and
$\hat{\Theta}_{A,A}$ and $\hat{\Omega}_{A,A}$ be estimators
of $\Theta_{A,A}$ and $\Omega_{A,A}$ defined in (\ref{estimator}).
Let $\delta\ge1$.
Suppose $s\leq c_{0}n/\log p$ for a sufficiently small constant $c_{0}>0$.
\begin{longlist}[(iii)]
\item[(i)] Suppose that conditions (\ref{cond-1}), (\ref{cond-2}) and
(\ref{cond-3}) hold
with $C_0$ and $\epsilon_\Omega$. Then
%
\begin{equation}
\max_{A:A=\{i,j\}}\mathbb{%
P} \biggl\{ \bigl\llVert \hat{
\Theta}_{A,A}-\Theta_{A,A}^{\mathrm{ora}}\bigr\rrVert
_{\infty}>C_{1}s\frac{\delta\log p}{n} \biggr\} \le6\epsilon
_\Omega p^{-\delta+1} + \frac{4p^{-\delta+1}}{\sqrt{2\log p}} \label{result1}
\end{equation}
with a positive constant $C_{1}$ depending on $\{C_0,\max_{m\in A=\{
i,j\}}\theta_{mm}\}$ only, and
%
\begin{equation}
\max_{A:A=\{i,j\}}\mathbb{%
P} \biggl\{ \bigl
\llVert \hat{\Omega}_{A,A}-\Omega_{A,A}^{\mathrm{ora}}\bigr
\rrVert _{\infty}>C_{1}^{\prime}s\frac{\delta\log p}{n} \biggr
\} \le6\epsilon_\Omega p^{-\delta+1} + \frac{4p^{-\delta+1}}{\sqrt
{2\log p}}
\label{result1a}
\end{equation}
%
with a constant $C_{1}^{\prime}>0$ depending on
$\{c_0C_1,\max_{m\in A=\{i,j\}} \{\omega_{mm},\theta_{mm}
\}\}
$ only.

\item[(ii)]
Let $\lambda= ( 1+\varepsilon ) \sqrt{%
\frac{2\delta\log p}{n}}$ with $\varepsilon>0$ in (\ref{Scaled Lasso}),
$\hat{\beta}_{A^c,A}$ be the scaled lasso estimator (\ref{Scaled Lasso})
or the LSE after the scaled lasso selection (\ref{Scaled Lasso LSE}).
Then (\ref{cond-1}), (\ref{cond-2}) and~(\ref{cond-3}), and thus
(\ref
{result1}) and (\ref{result1a}),
hold for all $\Omega\in\mathcal{G}^{\ast}(M,s,\lambda)$ with
a certain constant $C_0$ depending on $\{\varepsilon,c_{0},M\}$ only and
%
\begin{equation}
\max_{\Omega\in\mathcal{G}^{\ast}(M,s,\lambda)}\epsilon_\Omega = o (1 ). \label{result1b}
\end{equation}

\item[(iii)] Let $\kappa_{ij}^{\mathrm{ora}}=\sqrt{n}(\omega
_{ij}^{\mathrm{ora}}-\omega_{ij})/
\sqrt{\omega_{ii}\omega_{jj}+\omega_{ij}^{2}}$.
There exist constants $D_{1}$ and $\vartheta\in (
0,\infty ) $, and four marginally standard normal random
variables $Z^{\prime}$,
$Z_{kl}$, where $kl=ii,ij,jj$, such that whenever $\llvert
Z_{kl}\rrvert \leq\vartheta\sqrt{n}$ for all $kl$, we have
%
\begin{equation}
\bigl\llvert \kappa_{ij}^{\mathrm{ora}}-Z^{\prime}\bigr\rrvert
\leq\frac
{D_{1}}{\sqrt{%
n}} \bigl( 1+Z_{ii}^{2}+Z_{ij}^{2}+Z_{jj}^{2}
\bigr). \label{TailOracle2}
\end{equation}
Moreover, $Z^{\prime}$ can be defined
as a linear combination of $Z_{kl}$, $kl=ii,ij,jj$.
\end{longlist}
\end{theorem}

Theorem~\ref{Main} immediately yields the following results of estimation
and inference for $\omega_{ij}$.

\begin{theorem}
\label{Main1} Let $\hat{\Omega}_{A,A}$ be the estimator of $\Omega_{A,A}$
in (\ref{estimator}) with the components of $\hat{\bolds \epsilon}_A$
being the
estimated residuals (\ref{estimator1}) of (\ref{Scaled Lasso})
or (\ref{Scaled Lasso LSE}). Set $\lambda= ( 1+\varepsilon
 )
\sqrt{\frac{2\delta\log p}{n}}$ in (\ref{Scaled Lasso}) with certain
$\delta\geq1$ and $\varepsilon>0$.
Suppose $s\le c_0 n/\log p$ for a
sufficiently small constant $c_0>0$.
For any small constant $\epsilon_0 >0$, there exists a constant $%
C_{2}=C_{2} ( \epsilon_0,\varepsilon,c_0,M ) $ such that
%
\begin{equation}
\max_{\Omega\in\mathcal{G}^*(M,s,\lambda)} \max_{1\le i\le j\le p}
\mathbb{P}%
 \biggl\{ \llvert \hat{\omega}_{ij}-
\omega_{ij}\rrvert >C_{2}\max \biggl\{ s\frac{\log p}{n},
\sqrt{\frac{1}{n}} \biggr\} \biggr\} \leq \epsilon_0 .
\label{resultomega}
\end{equation}
Moreover, there exists a constant $C_{3}=C_3 (\delta,\varepsilon
,c_0,M )$ such that
%
\begin{equation}
\max_{\Omega\in\mathcal{G}^*(M,s,\lambda)} \mathbb{P} \biggl\{ \llVert \hat{%
\Omega}-\Omega\rrVert _{\infty}>C_{3}\max \biggl\{ s
\frac{\log
p}{n},%
\sqrt{\frac{\log p}{n}} \biggr\} \biggr\} = o
\bigl( p^{-\delta
+3} \bigr) . \label{result3a}
\end{equation}
Furthermore, $\hat{\omega}_{ij}$ is asymptotically efficient with a
consistent variance estimate
%
\begin{equation}
\sqrt{nF_{ij}} ( \hat{\omega}_{ij}-\omega_{ij} )
\overset{D} {%
\rightarrow}\mathcal{N} ( 0,1 ), \qquad {\hat
F}_{ij}/F_{ij}\to 1, \label{result 4a}
\end{equation}
uniformly for all $i,j$ and $\Omega\in\mathcal{G}^*(M,s,\lambda)$,
provided that $s=o (\sqrt{n}/\log
p )$, where
\[
F_{ij}=\bigl(\omega_{ii}\omega_{jj}+
\omega_{ij}^{2}\bigr)^{-1},\qquad {\hat F}_{ij}=
\bigl({\hat\omega}_{ii}{\hat\omega}_{jj}+{\hat\omega
}_{ij}^{2}\bigr)^{-1}.
\]
\end{theorem}

\begin{remark}
The upper bounds $\max \{ s\frac{\log p}{n},\sqrt{\frac
{1}{n}} \} $
and $\max \{ s\frac{\log p}{n},\sqrt{\frac{\log p}{n}}
\} $ in
equations (\ref{resultomega}) and (\ref{result3a}), respectively, are
shown to
be rate-optimal in Section~\ref{low.sec}.
\end{remark}

\begin{remark}
\label{lambdachoice}The choice of $\lambda= ( 1+\varepsilon
 )
\sqrt{\frac{2\delta\log p}{n}}$ is common in the literature, but can
be too
big and too conservative, which usually leads to some estimation bias
in practice.
Let $L_n(t)$ be the negative quantile function of $\mathcal{N}(0,1/n)$,
which satisfies
$L_n(t)\approx\sqrt{(2/n)\log p}$.
In Sections~\ref{regression-revisit} and \ref{Discu} we show the value
of $\lambda$ can be reduced to
$ ( 1+\varepsilon )L_n(k/p)$ when $\delta\vee k=o(\sqrt
{n}/\log p)$.
\end{remark}

\begin{remark}
\label{Partial Correlation}In Theorems \ref{Main} and \ref{Main1},
our goal
is to estimate each entry $\omega_{ij}$ of the precision matrix
$\Omega$.
Sometimes it is more natural to consider estimating the partial
correlation $%
r_{ij}=-\omega_{ij}/(\omega_{ii}\omega_{jj})^{1/2}$ between $Z_{i}$
and $%
Z_{j}$. Let $\hat{\Omega}_{A,A}$ be estimator of $\Omega_{A,A}$
defined in (%
\ref{estimator}). Our estimator of partial correlation $r_{ij}$ is defined
as $\hat{r}_{ij}=-\hat{\omega}_{ij}/(\hat{\omega}_{ii}\hat{\omega
}%
_{jj})^{1/2}$. Then the results above can be easily extended to the
case of
estimating $r_{ij}$. In particular, under the assumptions of Theorem~\ref{Main1}, the estimator $\hat{r}_{ij}$ is asymptotically
efficient: $
\sqrt{n(1-r_{ij}^{2})^{-2}}(\hat{r}_{ij}-r_{ij})$ converges to
$\mathcal
{N}%
 ( 0,1 ) $ when $s=o (\sqrt{n}/\log p )$. This
asymptotic normality result was stated
as Corollary $1$ in \citet{SunZhang2012SinicaComment} without proof.
\end{remark}

The following theorem extends Theorems \ref{Main} and \ref{Main1} to
the estimation of $\zeta ( \Omega_{B,B}^{-1} ) $, a smooth
functional
of $\Omega_{B,B}^{-1}$ for a fixed size subset $B$.
Assume that $\zeta\dvtx %
\mathbb{R}^{\llvert  B\rrvert \times\llvert  B\rrvert
}\rightarrow\mathbb{R}$ is a unit Lipschitz function in a neighborhood
$%
 \{ G\dvtx |\!|\!|G-\Omega_{B,B}^{-1}|\!|\!|\leq\kappa \} $, that is,
%
\begin{equation}
\bigl\llvert \zeta ( G ) -\zeta \bigl( \Omega _{B,B}^{-1}
\bigr) \bigr\rrvert \leq\bigl|\!\bigl|\!\bigl|G-\Omega_{B,B}^{-1}\bigr|\!\bigr|\!\bigr|.
\label{lips}
\end{equation}

\begin{theorem}
\label{Main 2} Let $\hat{\zeta}$ be the estimator of $\zeta$ defined
in (\ref%
{estimator3}) with the components of $\hat{\bolds {\epsilon}}_{B}$
being the
estimated residuals (\ref{estimator1}) of the estimators (\ref{Scaled Lasso})
or (\ref{Scaled Lasso LSE}). Set the penalty level $\lambda= (
1+\varepsilon )
\sqrt{\frac{2\delta\log p}{n}}$ in (\ref{Scaled Lasso}) with certain
$\delta\geq1$ and $\varepsilon>0$.
Suppose $s\le c_0 n/\log p$ for a
sufficiently small constant $c_0>0$. Then
\begin{equation}
\max_{\Omega\in\mathcal{G}^*(M,s,\lambda)} \mathbb{P} \biggl\{ \bigl\llvert \hat{%
\zeta}-\zeta^{\mathrm{ora}}\bigr\rrvert >C_{1}s\frac{\log p}{n}
\biggr\} = o \bigl( \llvert B\rrvert p^{-\delta+1} \bigr) , \label{general result}
\end{equation}
with a constant $C_{1}=C_1(\varepsilon,c_0,M,\llvert  B\rrvert )$.
Furthermore, $\hat{\zeta}$ is asymptotically efficient
%
\begin{equation}
\sqrt{nF_{\zeta}} ( \hat{\zeta}-\zeta ) \overset {D} {
\rightarrow}%
\mathcal{N} ( 0,1 ) , \label{general result 2}
\end{equation}
when $\Omega\in\mathcal{G}^*(M,s,\lambda)$ and $s=o (\sqrt
{n}/\log
p )$, where $F_{\zeta}$ is the Fisher information of estimating
$\zeta
$ for the Gaussian model $\mathcal{N} ( 0,\Omega
_{B,B}^{-1} ) $.
\end{theorem}

The results in this section can be easily
extended to the weak $\ell_{q}$ ball with $0 < q<1$ to model the
sparsity of
the precision matrix. A weak $\ell_{q}$ ball of radius $c$ in $\mathbb{R}^{p}$
is defined as follows:
\[
B_{q} ( c ) = \bigl\{ \xi\in\mathbb{R}^{p}\dvtx \llvert
\xi _{ (
j ) }\rrvert ^{q} \leq cj^{-1},\mbox{ for all
}j=1,\ldots ,p \bigr\},
\]
where $\llvert \xi_{ ( 1 ) }\rrvert \geq\llvert \xi
_{ ( 2 ) }\rrvert \geq\cdots\geq\llvert \xi
_{ (
p ) }\rrvert $. Let
%
\begin{equation}
\mathcal{G}_{q}(M,k_{n,p})=\lleft\{ %
\begin{array} {c}\displaystyle \Omega= ( \omega_{ij} ) _{1\leq i,j\leq p}\dvtx
\omega_{\cdot
j}\in B_{q} ( k_{n,p} ) ,
\\[2pt]
\mbox{and }\displaystyle 1/M\leq\lambda_{\min} ( \Omega ) \leq \lambda _{\max
}
( \Omega ) \leq M%
\end{array} %
 \rright\}.
\label{weakLqBall}
\end{equation}
Since $\xi\in B_q(k)$ implies $\sum_j\min \{1,|\xi_j|/\lambda
\}\le\lfloor
k/\lambda^q\rfloor+ \{q/(1-q)\}k^{1/q} \lfloor
k/\break \lambda^q\rfloor^{1-1/q}/\lambda$,
%
\begin{equation}
\label{LqCappedL1} \mathcal{G}_{q}(M,k_{n,p})\subseteq
\mathcal{G}^*(M,s,\lambda),\qquad 0\le q < 1,
\end{equation}
when $C_{q}k_{n,p}/\lambda^q\le s$, where $C_q=1+q 2^{1/q-1}/(1-q)$ for
$%
0<q<1$ and $C_0=1$. We state the extension in the following corollary.

\begin{corollary}\label{UpperBdextenToLqball}
The conclusions of Theorems \ref{Main}, \ref{Main1}
and \ref{Main 2} hold with $\mathcal{G}^*(M, \break s,\lambda)$ replaced by $
\mathcal{G}_{q}(M,k_{n,p})$ and $s$ by $k_{n,p}(n/\log p)^{q/2}$,
$0\le q<1$.
\end{corollary}
%

\subsection{Lower bound}
\label{low.sec}

In this section, we derive a lower bound for estimating $\omega_{ij}$ over
the matrix class $\mathcal{G}_{0}(M,k_{n,p})$ defined in (\ref%
{sparseparaspace}). Assume that
%
\begin{equation}
p\geq k_{n,p}^{\nu}\qquad\mbox{with }\nu>2 \label{techassump2}\vadjust{\goodbreak}
\end{equation}
and
%
\begin{equation}
3 \leq k_{n,p}\leq C_{0}\frac{n}{\log p} \label{techassump}
\end{equation}
for some $C_{0}>0$. Theorem~\ref{sparseMinimaxlwSupport} below implies that
the assumption $k_{n,p}\frac{\log p}{n}\rightarrow0$ is necessary for
consistent estimation of any single entry of $\Omega$.

We carefully construct a finite collection of distributions
$\mathcal{G}_{0}\subset\mathcal{G}_{0}(M,k_{n,p})$
and apply Le Cam's method to show
that for any estimator $\hat{\omega}_{ij}$,
%
\begin{equation}
\sup_{\mathcal{G}_{0}}\mathbb{P} \biggl\{ \llvert \hat{\omega
}_{ij}-\omega _{ij}\rrvert >C_{1}k_{n,p}
\frac{\log p}{n} \biggr\} \mathbb {\rightarrow}%
1, \label{LW2}
\end{equation}
for some constant $C_{1}>0$. It is relatively easy to establish the
parametric lower bound $\sqrt{\frac{1}{n}}$. These two lower bounds together
immediately yield Theorem~\ref{sparseMinimaxlwSupport} below.

\begin{theorem}
\label{sparseMinimaxlwSupport} Suppose we observe independent and
identically distributed $p$-variate Gaussian random variables $%
X^{(1)},X^{(2)},\ldots,X^{(n)}$ with zero mean and precision matrix
$\Omega
= ( \omega_{kl} ) _{p\times p}\in\mathcal{G}_{0}(M,k_{n,p})$.
Under assumptions (\ref{techassump2}) and (\ref{techassump}), we have the
following minimax lower bounds:
%
\begin{equation}
\inf_{\hat{\omega}_{ij}}\sup_{\mathcal{G}_{0}(M,k_{n,p})}\mathbb {P} \biggl\{
\llvert \hat{\omega}_{ij}-\omega_{ij}\rrvert >\max \biggl\{
C_{1} 
\frac{k_{n,p}\log p}{n},C_{2}\sqrt{
\frac{1}{n}} \biggr\} \biggr\} >c_{1}>0%
\label{sparserateSupport}
\end{equation}
and
%
\begin{equation}\qquad
\inf_{\hat{\Omega}}\sup_{\mathcal{G}_{0}(M,k_{n,p})}\mathbb {P} \biggl\{
\llVert \hat{\Omega}-\Omega\rrVert _{\infty}>\max \biggl\{
C_{1}^{\prime}\frac{k_{n,p}\log p}{n},C_{2}^{\prime}
\sqrt{\frac
{\log
p}{n}}%
 \biggr\} \biggr\} >c_{2}>0,
\label{sparserateSupport1}
\end{equation}
where $c_1,c_2,C_{1}$, $C_{2}$, $C_{1}^{\prime}$ and $C_{2}^{\prime}$
are positive
constants depending on $M$, $\nu$ and $C_{0}$ only.
\end{theorem}

\begin{remark}
\label{assumpCnp} The lower bound $\frac{k_{n,p}\log p}{n}$
in Theorem~\ref{sparseMinimaxlwSupport} shows that estimation of sparse
precision matrix can be very different from estimation of sparse covariance
matrix. The sample covariance always gives a parametric rate of estimation
for every entry $\sigma_{ij}$. But for estimation of sparse
precision matrix, when $k_{n,p}\gg\frac{\sqrt{n}}{\log p%
}$, Theorem~\ref{sparseMinimaxlwSupport} implies that it is
impossible to obtain the parametric rate.
\end{remark}

\begin{remark}
\label{extenToLqball} Since $\mathcal{G}_{0}(M,k_{n,p})\subseteq
\mathcal{G}%
^*(M,k_{n,p},\lambda)$ by the definitions in (\ref{sparseparaspace})
and~(%
\ref{StarBall}), Theorem~\ref{sparseMinimaxlwSupport} also provides the
lower bound for the larger class. Similarly, Theorem~\ref%
{sparseMinimaxlwSupport} can be easily extended to the weak $\ell_{q}$
ball, $%
0<q<1$, defined in (\ref{weakLqBall}) and the capped-$\ell_1$ ball defined
in (\ref{StarBall}). 
For these parameter spaces, in the proof of Theorem~\ref%
{sparseMinimaxlwSupport} we only need to define $\mathcal{H}$ as the
collection of all $p\times p$ symmetric matrices with exactly $ (
k_{n,p} ( \frac{n}{\log p} ) ^{q/2}-1 ) $ rather than
$ (
k_{n,p}-1 ) $ elements equal to $1$ between the third and the last
elements on the first row (column) and the rest all zeros. Then it is easy
to check that the sub-parameter space $\mathcal{G}_{0}$ in (\ref
{G0}) is
indeed in $\mathcal{G}_{q}(M,k_{n,p})$. Now under assumptions $p\geq
 (
k_{n,p} ( \frac{n}{\log p} ) ^{q/2} ) ^{v}$ with $\nu>2$
and $%
k_{n,p}\leq C_{0} ( \frac{n}{\log p} ) ^{1-q/2}$, we have the
following minimax lower bounds:
\[
\inf_{\hat{\omega}_{ij}}\sup_{\mathcal{G}_{q}(M,k_{n,p})}\mathbb {P} \biggl\{
\llvert \hat{\omega}_{ij}-\omega_{ij}\rrvert >\max \biggl\{
C_{1}k_{n,p} \biggl( \frac{\log p}{n} \biggr)
^{1-q/2},C_{2}\sqrt {\frac
{1}{n}}%
 \biggr\}
\biggr\} >c_{1}>0
\]
and
\[
\inf_{\hat{\Omega}}\sup_{\mathcal{G}_{q}(M,k_{n,p})}\mathbb {P} \biggl\{
\llVert \hat{\Omega}-\Omega\rrVert _{\infty}>\max \biggl\{
C_{1}^{\prime}k_{n,p} \biggl( \frac{\log p}{n}
\biggr) ^{1-q/2},C_{2}^{\prime}%
\sqrt{
\frac{\log p}{n}} \biggr\} \biggr\} >c_{2}>0.
\]
These lower bounds match the upper bounds in
Corollary~\ref{UpperBdextenToLqball} for the proposed estimator.
\end{remark}

\section{Applications}
\label{app}

The asymptotic normality result is applied to obtain rate-optimal estimation
of the precision matrix under various matrix $\ell_{w}$ norms, to
recover the
support of $\Omega$ adaptively and to estimate latent graphical models
without the need of the irrepresentability condition or the $\ell_{1}$
constraint
of the precision matrix commonly required in literature. In our
procedure, we
first obtain an asymptotically normal estimation and then thresholding.
We thus call it ANT.

\subsection{ANT for adaptive support recovery}

The support recovery of precision matrix has been studied by several papers.
See, for example, \citet{FHT08}, \citet{d'ABE08}, \citet{RBLZ08}, \citet
{RWRY08}, \citet{CLL11} and \citet{CLZ12}. In these works, the theoretical properties
of the graphical lasso (GLasso), CLIME and ACLIME on the support recovery
were obtained. \citet{RWRY08} studied the theoretical properties of GLasso,
and showed that GLasso can correctly recover the support under
a strong irrepresentability condition and a uniform signal strength
condition $\min_{(i,j)\dvtx \omega_{ij}\neq0}|\omega_{ij}|\geq c\sqrt
{\frac{\log p}{n}}$ for some $c>0$. \citet%
{CLL11} do not require irrepresentability conditions, but need to
assume that
$\min_{(i,j)\dvtx \omega_{ij}\neq0}|\omega_{ij}|\geq CM_{n,p}^{2}\sqrt{
\frac{\log p}{n}}$, where $M_{n,p}$ is the matrix $\ell_{1}$ norm of
$\Omega$.
In \citet{CLZ12}, they weakened the condition to $\min_{(i,j)\dvtx \omega
_{ij}\neq0}|\omega_{ij}|\geq CM_{n,p}\sqrt{\frac{\log p}{n}}$, but the
threshold level there is $\frac{C}{2}M_{n,p}\sqrt{\frac{\log p}{n}}$,
where $%
C$ is unknown and $M_{n,p}$ can be very large, which makes the support
recovery procedure there impractical.

In this section we introduce an adaptive support recovery procedure
based on
the variance of the oracle estimator of each entry $\omega_{ij}$ to recover
the sign of nonzero entries of $\Omega$ with high probability. The lower
bound condition for $\min_{(i,j)\dvtx \omega_{ij}\neq0}|\omega_{ij}|$ is
significantly weakened. In particular, we remove the unpleasant matrix $
\ell_{1} $ norm $M_{n,p}$. In Theorem~\ref{Main1}, when the
precision matrix
is sparse enough $s=o ( \frac{\sqrt{n}}{\log p} ) $, we
have the
asymptotic normality result for each entry $\omega_{ij}$, $i\neq j$,
that is,
\[
\sqrt{nF_{ij}} ( \hat{\omega}_{ij}-\omega_{ij} )
\overset{D} {%
\rightarrow}\mathcal{N} ( 0,1 ) ,
\]
where $F_{ij}= ( \omega_{ii}\omega_{jj}+\omega_{ij}^{2} ) ^{-1}$
is the Fisher information of estimating $\omega_{ij}$. The total
number of
edges is $p ( p-1 ) /2$. We may apply thresholding to $\hat
{\omega}%
_{ij}$ to correctly distinguish zero and nonzero entries. However, the
variance $\omega_{ii}\omega_{jj}+\omega_{ij}^{2}$ needs to be estimated.
We define the adaptive support recovery procedure as follows:
%
\begin{equation}
\label{adaptive thresholding} \hat{\Omega}_{\mathrm{thr}}=\bigl(\hat{\omega}_{ij}^{\mathrm{thr}}
\bigr)_{p\times p},
\end{equation}
where $\hat{\omega}_{ii}^{\mathrm{thr}}=\hat{\omega}_{ii}$ and $\hat{\omega
}%
_{ij}^{\mathrm{thr}}=\hat{\omega}_{ij}1\{|\hat{\omega}_{ij}|\geq{\hat{\tau
}}_{ij}\}$
for $i\neq j$ with
%
\begin{equation}
\label{thresholding level} {\hat{\tau}}_{ij}=\sqrt{\frac{2\xi_{0} ( \hat{\omega
}_{ii}\hat
{\omega}%
_{jj}+\hat{\omega}_{ij}^{2} ) \log p}{n}}.
\end{equation}
Here $\hat{\omega}_{ii}\hat{\omega}_{jj}+\hat{\omega}_{ij}^{2}$
is the
natural estimate of the asymptotic variance of $\hat{\omega}_{ij}$ defined
in (\ref{estimator}), and $\xi_{0}$ is a tuning parameter which can be
taken as fixed at any $\xi_{0}>2$. This thresholding estimator is adaptive.
The sufficient conditions in Theorem~\ref{support recovery} below for
support recovery are much weaker than other results in literature.

Define a thresholded population precision matrix as
%
\begin{equation}
\label{thresholded-Omega} \Omega_{\mathrm{thr}} =\bigl({\omega}_{ij}^{\mathrm{thr}}
\bigr)_{p\times p},
\end{equation}
where $\omega_{ii}^{\mathrm{thr}}=\omega_{ii}$ and $\omega_{ij}^{\mathrm{thr}}
=\omega_{ij}1 \{|\omega_{ij}|\geq\sqrt{8\xi( \omega
_{ii}\omega
_{jj}+\omega_{ij}^{2})(\log p)/n} \}$, with a certain $\xi>\xi_0$.
Recall that $E = E(\Omega) = \{(i,j)\dvtx \omega_{ij}\neq0\}$ is the
edge set
of the Gauss--Markov graph associated with the precision matrix $\Omega$.
Since $\Omega_{\mathrm{thr}}$ is composed of relatively large components of
$\Omega$,
$(V,E(\Omega_{\mathrm{thr}}))$ can be viewed as a graph of strong edges. Define
\[
\mathcal{S}(\Omega)=\bigl\{\operatorname{sgn}(\omega_{ij}),1\leq i,j\leq p\bigr\}.
\]
The following theorem shows that with high probability, ANT recovers
all the
strong edges without false recovery. Moreover, under the uniform signal
strength condition,
%
\begin{equation}
\llvert \omega_{ij}\rrvert \geq2\sqrt{\frac{2\xi (
\omega
_{ii}\omega_{jj}+\omega_{ij}^{2} ) \log p}{n}}\qquad
\forall \omega_{ij} \neq0; \label{lower bound for each entry}
\end{equation}
that is, $\Omega_{\mathrm{thr}}=\Omega$, and the ANT also recovers the sign
matrix $\mathcal{S%
}(\Omega)$.

\begin{theorem}
\label{support recovery} Let $\lambda= ( 1+\varepsilon )
\sqrt{%
\frac{2\delta\log p}{n}}$ for any $\delta\geq3$ and $\varepsilon
>0$. Let
$\hat{\Omega}_{\mathrm{thr}}$ be the ANT estimator defined in (\ref{adaptive
thresholding}) with $\xi_{0}>2$ in the thresholding level (\ref%
{thresholding level}). Suppose $\Omega\in\mathcal{G}^{\ast
}(M,s,\lambda)$
with $s=o ( \sqrt{n/\log p} ) $. Then
%
\begin{equation}
\lim_{n\rightarrow\infty}\mathbb{P} \bigl( E(\Omega _{\mathrm{thr}})
\subseteq E(\hat{%
\Omega}_{\mathrm{thr}})\subseteq E(\Omega) \bigr) =1,
\label{result: support recovery}
\end{equation}
where $\Omega_{\mathrm{thr}}$ is defined in (\ref{thresholded-Omega}) with
$\xi
>\xi
_{0}$. If in addition (\ref{lower bound for each entry}) holds, then
\begin{equation}
\lim_{n\rightarrow\infty}\mathbb{P} \bigl( \mathcal{S}(\hat {\Omega
}_{\mathrm{thr}})=%
\mathcal{S}(\Omega) \bigr) =1. \label{result: support recovery2}
\end{equation}
\end{theorem}

\subsection{ANT for adaptive estimation under the matrix \texorpdfstring{$\ell_{w}$}{ellw} norm}
In this section, we consider the rate of convergence under the matrix
$\ell_{w}$ norm.
To control the impact of extremely small tail probability of near
singularity of the low-dimensional
estimator $\hat{\Theta}_{A,A}$, we define a truncated version of the estimator
$\hat{\Omega}_{\mathrm{thr}}$ defined in (\ref{adaptive thresholding}),
%
\begin{equation}
\breve{\Omega}_{\mathrm{thr}}= \biggl(\hat{\omega}_{ij}^{\mathrm{thr}}
\min \biggl\{ 1,\frac
{\log p}{\llvert \hat{\omega}_{ij}\rrvert } \biggr\} \biggr)_{p\times p}.
\label{adaptive threshold- upper bounded}
\end{equation}
Theorem~\ref{matrixnorm} below follows mainly from the
convergence rate
under element-wise norm and the fact that the upper bound holds for the
matrix $%
\ell_{1}$ norm. This argument uses the inequality $|\!|\!|M|\!|\!|_{w}\leq
|\!|\!|M|\!|\!|_{1}$ for symmetric matrices $M$ and $1\leq w\leq\infty$, which
follows from
the Riesz--Thorin interpolation theorem; see, for example, \citet
{thorin1948convexity}.
Note that under the assumption $%
s^{2}=O ( n/\log p ) $, it can be seen from the equations (%
\ref{result1a}) and (\ref{TailOracle2}) in Theorem~\ref{Main} that
with high
probability the $\llVert \hat{\Omega}-\Omega\rrVert
_{\infty}$ is
dominated by $\llVert \Omega^{\mathrm{ora}}-\Omega\rrVert _{\infty
}=O_{p} ( \sqrt{\frac{\log p}{n}} ) $. From there the
details of
the proof is similar in nature to those of Theorem~3 in \citet
{CZh12Sparse} and
thus will be omitted due to the limit of space.

\begin{theorem}
\label{matrixnorm}Under the assumptions $s^2=O ( n/\log p ) $
and $%
n=O ( p^{\xi} ) $ with some $\xi>0$, the $\breve{\Omega}_{\mathrm{thr}}$
defined in (\ref{adaptive threshold- upper bounded}) with $\lambda
= (
1+\varepsilon ) \sqrt{\frac{2\delta\log p}{n}}$ for sufficiently
large $\delta\geq3+$ $2\xi$ and $\varepsilon>0$ satisfies, for all $
1\leq w\leq\infty$ and $k_{n,p}\le s$,
%
\begin{equation}
\sup_{\mathcal{G}_{0}(M,k_{n,p})}\mathbb{E}|\!|\!|\breve{\Omega}%
_{\mathrm{thr}}-
\Omega|\!|\!|_{w}^{2} \le\sup_{\mathcal{G}^*(M,k_{n,p},\lambda
)}
\mathbb{E%
}|\!|\!|\breve{\Omega}_{\mathrm{thr}}-\Omega|\!|
\!|_{w}^{2} \leq C s^{2}\frac{\log p}{n}.
\label{l0rate}
\end{equation}
\end{theorem}

\begin{remark}
It follows from equation (\ref{LqCappedL1}) that result (\ref
{l0rate}) also
holds for the classes of weak $\ell_p$ balls $\mathcal{G}_{q}(M,k_{n,p})$
defined in equation (\ref{weakLqBall}), with $s=C_q k_{n,p} (
\frac
{n}{%
\log p} ) ^{q/2}$,
%
\begin{equation}
\sup_{\mathcal{G}_{q}(M,k_{n,p})}\mathbb{E}|\!|\!|\breve{\Omega }_{\mathrm{thr}}-
\Omega |\!|\!|_{w}^{2}\leq Ck_{n,p}^{2}
\biggl( \frac{\log p}{n} \biggr) ^{1-q}. \label{lwrate}
\end{equation}
\end{remark}

\begin{remark}
\citet{CLZ12} showed that the rates obtained in equations (\ref{l0rate})
and (\ref%
{lwrate}) are optimal when $p\geq cn^{\alpha_0 }$ for some $\alpha_0
>1$ and $%
k_{n,p}=o ( n^{1/2} ( \log p ) ^{-3/2} ) $.
\end{remark}

\begin{remark}
Although the estimator $\breve{\Omega}_{\mathrm{thr}}$ is symmetric, it is not
guaranteed to be positive definite. It follows from Theorem~\ref{matrixnorm}
that $\breve{\Omega}_{\mathrm{thr}}$ is positive definite with high
probability. When
it is not positive definite, we can always pick the smallest $c_{a}\geq
0$ such that $c_{a}I+\breve{\Omega%
}_{\mathrm{thr}}$ is positive semidefinite. It is trivial to see that $%
(c_{a}+1/n)I+\breve{\Omega}_{\mathrm{thr}}$ is positive definite,
sparse and enjoys the same rate of convergence as $\breve{\Omega}_{\mathrm{thr}}$
for the loss functions considered in this paper.
\end{remark}

\subsection{Estimation and inference for latent variable graphical model}
\label{latent graphical model}

Chan\-drasekaran, Parrilo and Willsky (\citeyear{CPW12}) first proposed a very natural penalized estimation
approach and
studied its theoretical properties. Their work has been discussed and
appreciated
by several researchers, but it has never been clear if the conditions
in their
paper are necessary and the results optimal. \citet{RZ12} observed
that the support recovery boundary can be significantly improved from an
order of $\sqrt{\frac{p}{n}}$ to $\sqrt{\frac{\log p}{n}}$ under certain
conditions including a bounded $\ell_{1}$ norm constraint for the precision
matrix. In this section we extend the methodology and results in
Section~\ref%
{methodinf} to study latent variable graphical models. The results in
\citet%
{RZ12} are significantly improved under weaker assumptions.

Let $O$ and $H$ be two subsets of $ \{ 1,2,\ldots,p+h \} $ with
$\operatorname{Card} ( O ) =p$,\break  $\operatorname{Card} ( H ) =h$
and $%
O\cup H= \{ 1,2,\ldots,p+h \} $. Assume that $ (
X_{O}^{ (
i ) },X_{H}^{ ( i ) } ) $, $i=1,\ldots,n$, are
i.i.d. $%
 ( p+h ) $-variate Gaussian random vectors with a positive
covariance matrix $\Sigma_{ ( p+h ) \times (
p+h ) }$.
Denote the corresponding precision matrix by $\bar{\Omega}_{ (
p+h ) \times ( p+h ) }=\Sigma_{ ( p+h )
\times
 ( p+h ) }^{-1}$. We only have access to $ \{
X_{O}^{ (
1 ) },X_{O}^{ ( 2 ) },\ldots,X_{O}^{ ( n )
} \}
$, while $ \{ X_{H}^{ ( 1 ) },X_{H}^{ ( 2 )
}, \ldots
,X_{H}^{ ( n ) } \} $ are hidden and the number of latent
components is unknown. Write $\Sigma_{ ( p+h ) \times (
p+h ) }$ and $\bar{\Omega}_{ ( p+h ) \times (
p+h )
} $ as follows:
\[
\Sigma_{ ( p+h ) \times ( p+h ) }=\pmatrix{
\Sigma_{O,O} & \Sigma_{O,H}
\vspace*{2pt}\cr
\Sigma_{H,O} & \Sigma_{H,H}%
} \quad \mbox{and}\quad \bar{\Omega}_{ ( p+h ) \times
 (
p+h ) }=\pmatrix{ \bar{\Omega}_{O,O} & \bar{
\Omega}_{O,H}
\vspace*{2pt}\cr
\bar{\Omega}_{H,O} & \bar{\Omega}_{H,H}} ,
\]
where $\Sigma_{O,O}$ and $\Sigma_{H,H}$ are covariance matrices of $%
X_{O}^{ ( i ) }$ and $X_{H}^{ ( i ) }$,
respectively, and
from the Schur complement we have
%
\begin{equation}
\Sigma_{O,O}^{-1}=\bar{\Omega}_{O,O}-\bar{
\Omega}_{O,H}\bar {\Omega}%
_{H,H}^{-1}\bar{
\Omega}_{H,O}; \label{precision decomp}
\end{equation}
see, for example, \citet{MatAna}. Define
\[
S=\bar{\Omega}_{O,O},\qquad L=\bar{\Omega}_{O,H}\bar{\Omega
}_{H,H}^{-1}%
\bar{\Omega}_{H,O},
\]
and $h^{\prime}=\mathrm{rank} ( L )$. We note that
$h^{\prime}=\mathrm{rank} ( \bar{%
\Omega}_{O,H} ) \leq h$.

We focus on the estimation of $\Sigma_{O,O}^{-1}$ and $S$, as the
estimation of $L$
can be naturally carried out based on our results as in \citet{CPW12}
and \citet{RZ12}.
To make the problem identifiable we assume that $S$ is sparse, and the
observed and latent variables are weakly correlated in the following sense:
\begin{equation}
S= ( s_{ij} ) _{1\leq i,j\leq p},\qquad \max_{1\leq
j\leq
p}\sum
_{i=1}^{p}1 \{ s_{ij}\neq0 \} \leq
k_{n,p}, \label{Sstar}
\end{equation}
and that for certain $a_n\to0$
%
\begin{equation}
L= ( l_{ij} ) _{1\leq i,j\leq p},\qquad \max_j\sum
_i l_{ij}^2
\le(a_n/n)\log p. 
\label{elementwise bound}
\end{equation}
The sparseness of $S=\bar{\Omega}_{O,O}$ can be seen as inherited from
that of the full precision matrix $\bar{\Omega}_{ ( p+h )
\times
 ( p+h ) }$. It is particularly interesting for us to identify
the support of $S=\bar{\Omega}_{O,O}$ and make inference for each entry
of $S $.
The $\ell_2$ condition (\ref{elementwise bound}) on $L$ is of a
weaker form
than the $\ell_1$ and $\ell_\infty$ conditions imposed in \citet{RZ12}.
In addition, we assume that for some universal constant $M$,
%
\begin{equation}
1/M\leq\lambda_{\min}(\Sigma_{ ( p+h ) \times (
p+h )
})\leq\lambda_{\max}(
\Sigma_{ ( p+h ) \times (
p+h )
})\leq M, \label{eigenbd}
\end{equation}
%
which implies that both the covariance $\Sigma_{O,O}$ of
observations $X_{O}^{ ( i ) }$ and the sparse component
$S=\bar
{%
\Omega}_{O,O}$ have bounded spectrum. 

With a slight abuse of notation, we denote the precision matrix $\Sigma
_{O,O}^{-1}$ of $X_{O}^{ ( i ) }$ by $\Omega$ and its inverse
by $%
\Theta$. We propose the application of the methodology in Section~\ref{Methodology}
to i.i.d. observations $X^{ ( i ) }$ from $\mathcal{N}%
 ( 0,\Sigma_{O,O} ) $ with $\Omega= (
s_{ij}-l_{ij} )
_{1\leq i,j\leq p}$ by considering the following regression:
%
\begin{equation}
\mathbf{X}_{A}=\mathbf{X}_{O\setminus A}\bolds{\beta}+\bolds{
\epsilon}%
_{A} \label{regress}
\end{equation}
for $A= \{ i,j \} \subset O$ with $\bolds{\beta}=\Omega
_{O\setminus A,A}\Omega_{A,A}^{-1}$ and $\bolds{\epsilon
}_{A}\overset{%
\mathrm{i.i.d.}}{\sim}\mathcal{N} ( 0,\Omega_{A,A}^{-1} ) $.

To obtain the asymptotic normality result, condition (\ref{StarBall})
of Theorem~\ref{Main} requires
\[
\max_{j}\sum_{i=1}^{p}\min \biggl\{ 1,
\frac{\llvert
s_{ij}-l_{ij}\rrvert }{\lambda} \biggr\} =o \biggl( \frac{\sqrt
{n}}{\log p} \biggr) =o \biggl(
\frac{1}{\lam\sqrt{\log p}} \biggr)
\]
with $\lam\asymp\sqrt{(\log p)/n}$.
However, when $L$ is coherent [\citet{CR09}] in the sense of
$\{\max_{j}\sum_i |l_{ij}|\}^2 \asymp p \max_{i}\sum_jl_{ij}^2\asymp
p(a_n/n)\log p$,
\begin{eqnarray*}
\max_{j}\sum_{i=1}^{p}\min \biggl\{ 1,
\frac{\llvert
s_{ij}-l_{ij}\rrvert }{\lambda} \biggr\}& \ge&\max_{j}\sum_{s_{ij}=0}
\frac{|l_{ij}|}{\lam}\\
& \asymp&\frac{\sqrt{p}-\sqrt{k_{n,p}}}{\lam} \biggl(\frac{a_n\log
p}{n}
\biggr)^{1/2} \asymp\sqrt{a_np}.
\end{eqnarray*}
%
Thus the conditions of Theorem~\ref{Main} are not satisfied for the
latent variable
graphical model when $a_np(\log p)^{2}\ge n$. 
We overcome the difficulty through a new analysis.

\begin{theorem}\label{Latent Graphical Model Inference}
Let $\hat{\Omega}_{A,A}$ be the estimator of $\Omega_{A,A}$ defined in
(\ref{estimator})
with $A=\{i,j\}$ for the regression (\ref{regress}), where the components
of $\hat{\bolds \epsilon}_A$ are the estimated residuals of (\ref
{Scaled Lasso})
or (\ref{Scaled Lasso LSE}). Let $\lambda= ( 1+\varepsilon
 ) \sqrt{\frac{2\delta\log p}{n}}$ for certain $\delta\geq1$
and $%
\varepsilon>0$. Under assumptions (\ref{Sstar})--(\ref{eigenbd}) 
and $k_{n,p}\le c_0 n/\log p$ with a small $c_0$, we have
\[
\mathbb{P} \bigl\{ \llvert \hat{\omega}_{ij}-\omega_{ij}
\rrvert >C_{3}\max \bigl\{ k_{n,p}n^{-1}\log
p,n^{-1/2} \bigr\} \bigr\} = o \bigl( p^{-\delta+1} \bigr)
\]
for a certain constant $C_3$, and
\[
\mathbb{P} \bigl\{ \llvert \hat{\omega}_{ij}-s_{ij}\rrvert
>C_{3} \max \bigl\{ k_{n,p}n^{-1}\log
p,n^{-1/2},\sqrt{(a_n/n)\log p} \bigr\} \bigr\} = o \bigl(
p^{-\delta+1} \bigr).
\]
If the condition on $k_{n,p}$ is strengthened to $k_{n,p}=o (
\frac
{\sqrt{n}}{\log p} )$, then
%
\begin{equation}
\sqrt{\frac{n}{\omega_{ii}\omega_{jj}+\omega_{ij}^{2}}} ( \hat {\omega}%
_{ij}-
\omega_{ij} ) \overset{D} {\rightarrow}\mathcal{N} ( 0,1 ).
\label{result2a for Latent}
\end{equation}
\end{theorem}

\begin{remark}
If, in addition, $\ell_{ij}=o(n^{-1/2})$, then (\ref{result2a for
Latent}) implies
%
\begin{equation}
\sqrt{\frac{n}{\omega
_{ii}\omega_{jj}+\omega_{ij}^{2}}} ( \hat{\omega }_{ij}-s_{ij}
) \overset{D} {\rightarrow}\mathcal{N} ( 0,1 ). \label{result2 for Latent}
\end{equation}
\end{remark}

Define the adaptive thresholded estimator $\hat{\Omega}_{\mathrm{thr}}=(\hat
{\omega}%
_{ij}^{\mathrm{thr}})_{p\times p}$ as in (\ref{adaptive thresholding}) and
(\ref%
{thresholding level}). Following the proof of Theorems \ref{support recovery}
and \ref{matrixnorm}, we are able to obtain the following results. We shall
omit the proof due to the limit of space.

\begin{theorem}
\label{support recovery Latent} Let $\lambda= ( 1+\varepsilon
 )
\sqrt{\frac{2\delta\log p}{n}}$ for some $\delta\geq3$ and
$\varepsilon
>0 $ in (\ref{Scaled Lasso}). Assume assumptions (\ref{Sstar})--(\ref
{eigenbd})
hold. Then:
\begin{longlist}[(ii)]
\item[(i)] Under the assumptions $k_{n,p}=o ( \sqrt{\frac
{n}{\log
p}}%
 ) $ and
\[
\llvert s_{ij}\rrvert \geq2\sqrt{\frac{2\xi_{0} (
\omega
_{ii}\omega_{jj}+\omega_{ij}^{2} ) \log p}{n}}\qquad \forall
s_{ij} \neq0
\]
for some $\xi_{0}>2$, we have
\[
\lim_{n\rightarrow\infty}\mathbb{P} \bigl( \mathcal{S}(\hat {\Omega
}_{\mathrm{thr}})=%
\mathcal{S}(S) \bigr) =1.
\]

\item[(ii)] Under the assumption $k_{n,p}^{2}=O ( n/\log p
) $ and
$n=O ( p^{\xi} ) $ with some $\xi>0$, the $\breve{\Omega}_{\mathrm{thr}}$
defined in (\ref{adaptive threshold- upper bounded}) with sufficiently large
$\delta\geq3+$ $2\xi$ satisfies, for all $1\leq w\leq\infty$,
\[
\mathbb{E}|\!|\!|\breve{\Omega}_{\mathrm{thr}}-S|\!|\!|_{w}^{2}
\leq Ck_{n,p}^{2}\frac
{\log p}{%
n}.
\]
\end{longlist}
\end{theorem}

\section{Regression revisited}\label{regression-revisit}
The key element of our analysis is to establish (\ref{cond-1}), (\ref
{cond-2}) and (\ref{cond-3})
for the scaled lasso estimator (\ref{Scaled Lasso})
and the LSE after the scaled lasso selection (\ref{Scaled Lasso LSE}).
The existing literature has provided theorems and arguments to carry
out this task.
However, several issues still require extension of existing results or
explanation and modification
of existing proofs.
For example, the LSE after model selection is not as well understood as
the lasso,
and biased regression models are typically studied inexplicably, if at all.
Another issue is that the penalty level used in
theorems in previous sections could be too large for good numerical performance,
especially for $\delta\ge3$ in (\ref{result3a}) of Theorems \ref
{Main1} and
Theorems \ref{support recovery}, \ref{matrixnorm} and \ref{support
recovery Latent}.
These issues were addressed in previous versions of this paper
(\arxivurl{arXiv:1309.6024}) in separate lemmas.
In this section, we provide a streamlined presentation of these
regression results required in our analysis.

Let $\widetilde{\mathbf{ X}}=(\widetilde{\mathbf  X}_1,\ldots,\widetilde{\mathbf{X}}_{\widetilde p})$
be an $n\times{\widetilde p}$ standardized design matrix with $\|
\widetilde{\mathbf  X}_k\|^2=n$ for all
$k=1,\ldots,{\widetilde p}$, and $\widetilde{\mathbf  Y}$ be a response
vector satisfying
%
\begin{equation}
\label{LM-a} 
\widetilde{\mathbf Y} |\widetilde{\mathbf X}
\sim\mathcal{N} \bigl(\widetilde{\mathbf X}\bolds{\gamma}, \sigma ^2
\mathbf{I}_{n\times n} \bigr).
\end{equation}
%
For the scaled lasso $ \{\hat{\bolds {\beta}}_{m},\hat{\theta
}_{mm}^{1/2} \}$ in (\ref{Scaled Lasso}),
$ \{\Dbar_{A^c}^{1/2}\hat{\bolds {\beta}}_{m},\hat{\theta
}_{mm}^{1/2} \}$ can be written as
%
\begin{equation}
\{\hat{\bolds \gamma},\hat{\sigma} \}=\argmin_{\{
\bolds
{\gamma},\sigma\}} \biggl\{
\frac{ \Vert\widetilde{\mathbf  Y} -
\widetilde{\mathbf{X}}\bolds{\gamma}  \Vert^{2}}{2n\sigma} +
\frac{\sigma}{2} +\lambda_0 \llVert \bolds{
\gamma} \rrVert _{1} \biggr\}, \label{standard Scaled Lasso}
\end{equation}
with $m\in A=\{i,j\}$, $\widetilde{\mathbf  X}=\mathbf{X}_{A^c}\Dbar
_{A^c}^{-1/2}$,
$\Dbar= \operatorname{diag}(\mathbf{X}^T\mathbf{X}/n)$,
$\widetilde{\mathbf  Y} = \mathbf{X}_m$ and $\bolds{\gamma} = \Dbar
_{A^c}^{1/2}\bolds{\beta}_{m}$.
For the LSE after model selection in (\ref{Scaled Lasso LSE}),
$ \{\Dbar_{A^c}^{1/2}\hat{\bolds {\beta}}_{m},\hat{\theta
}_{mm}^{1/2} \}$ can be written as
%
\begin{equation}
\bigl\{\hat{\bolds \gamma}^{\mathrm{lse}},\hat{\sigma}^{\mathrm{lse}} \bigr\} =
\argmin _{\{\bolds{\gamma},\sigma\}} \biggl\{\frac{ \Vert\widetilde{\mathbf  Y} - \widetilde{\mathbf{X}}\bolds{\gamma}  \Vert^{2}}{2n\sigma} +\frac{\sigma}{2}\dvtx
\operatorname{ supp}(\bolds{\gamma}) \subseteq\operatorname{ supp}(\hat{\bolds
\gamma}) \biggr\}. \label{Scaled Lasso LSEa}
\end{equation}
%
Moreover, for both estimators, conditions (\ref{cond-1}), (\ref
{cond-2}) and (\ref{cond-3}) are consequences of
%
\begin{eqnarray}
\label{cond-1a} 
\mathbb{P} \bigl\{\bigl\llVert \widetilde{\mathbf X}\bigl(\hat{\bolds \gamma}-
\bolds{\gamma}^{\mathrm{target}}\bigr)\bigr\rrVert ^2 \le
C_0s\bigl(\sigma^{\mathrm{ora}}\bigr)^2\delta\log{
\widetilde p} \bigr\}&\ge& 1-{\widetilde p}^{1-\delta}{\widetilde
\epsilon}_0,
\\
\label{cond-2a} 
\mathbb{P} \bigl\{\bigl\llVert \hat{\bolds \gamma}- \bolds{\gamma
}^{\mathrm{target}}\bigr\rrVert _1 \le C_0s
\sigma^{\mathrm{ora}}\sqrt{\delta(\log{\widetilde p})/n} \bigr\} &\ge& 1- {
\widetilde p}^{1-\delta}{\widetilde\epsilon}_0
\end{eqnarray}
and
%
\begin{eqnarray}
\label{cond-3a} 
&&
\mathbb{P} \biggl\{ \biggl|\frac{\hat{\sigma}}{\sigma
^{\mathrm{ora}}}-1\biggr | \le C_0s\delta(\log{
\widetilde p})/n \le1/2 \biggr\}\ge1- {\widetilde p}^{1-\delta}{\widetilde
\epsilon}_0,
\end{eqnarray}
with $\sigma^{\mathrm{ora}} = \llVert \widetilde{\mathbf  Y}
- \widetilde{\mathbf  X}\bolds{\gamma}^{\mathrm{target}} \rrVert /\sqrt{n}$,
$\bolds{\gamma}^{\mathrm{target}}=\bolds{\gamma}$ and $\delta\ge1$, %
provided that $C_0$ is fixed and ${\widetilde\epsilon}_0\to0$
uniformly in $m\in A=\{i,j\}$ and
$\Omega$ in the class in (\ref{StarBall}); see Proposition~\ref{prop-probab}.
In the latent variable graphical model,
Theorems \ref{Latent Graphical Model Inference} and \ref{support
recovery Latent}
require (\ref{cond-1a}), (\ref{cond-2a}) and (\ref{cond-3a}) for a
certain sparse
$\bolds{\gamma}^{\mathrm{target}}$ in a biased linear model
when (\ref{LM-a}) does not provide a sufficiently sparse $\bolds
{\gamma}$.
We note that both
$\bolds{\gamma}$ and $\bolds{\gamma}^{\mathrm{target}}$ are allowed to be
random variables here.

To carry out an analysis of the lasso, one has to make a choice among different
ways of controlling the correlations between the design and noise
vectors in (\ref{LM-a}),
%
\begin{equation}
\label{v} \widetilde{\mathbf Z} = ({\widetilde Z}_1,\ldots,{
\widetilde Z}_{\widetilde p}) 
=
\frac{\widetilde{\mathbf  X}^T(\widetilde{\mathbf  Y} - \widetilde{\mathbf{X}}\bolds{\gamma})}{
\sqrt{n}\|\widetilde{\mathbf  Y} - \widetilde{\mathbf  X}\bolds
{\gamma
}\|}.
\end{equation}
%
A popular choice is to bound $\widetilde{\mathbf  Z}$ with the $\ell
_\infty$ norm as it is
the dual of the $\ell_1$ penalty. This has led to the sparse Riesz
[\citet{ZhangH08}], restricted eigenvalue [\citet{bickel2009simultaneous,
Koltchinskii09}],
compatibility [\citet{van2009conditions}],
cone invertibility [\citet{YZH10,ZhangZ12StatSc}] and other similar
conditions on the design matrix.
\citet{SZH12} took this approach to analyze
(\ref{standard Scaled Lasso}) and (\ref{Scaled Lasso LSEa})
with the compatibility and cone invertibility factors.
Another approach is to control the sparse $\ell_2$ norm of $\widetilde{\mathbf Z}$
to allow smaller penalty levels in the analysis; See \citet{Zhang09-l1}
and \citet{YZH10}
for analyses of the lasso and \citet{SZH122} for an analysis of
the scaled estimators (\ref{standard Scaled Lasso}) and (\ref{Scaled
Lasso LSEa}).

Here we take a different approach by using two threshold levels, a
smaller one to
bound an overwhelming majority of the components of
$\widetilde{\mathbf  Z}$ and a larger one to bound its $\ell_\infty$ norm.
This allows us to use both a small penalty level associated with the
smaller threshold
level and the compatibility condition.

For $\alpha\ge0$ and index sets $K$, the compatibility constant is
defined as
\[
\phi_{\mathrm{comp}} ( \alpha,K,\widetilde{\mathbf X} ) =\inf \biggl\{
\frac{\llvert  K\rrvert ^{1/2}\|\widetilde{\mathbf  X}u\|}{
n^{1/2}\llVert  u_{K}\rrVert _{1}}\dvtx u\in\mathcal{C} (\alpha,K ) ,u\neq0 \biggr\} ,
\]
where $\llvert  K\rrvert $ is the cardinality of $K$ and
$\mathcal{C} ( \alpha,K ) = \{ u\in\mathbb
{R}^{\widetilde p}\dvtx \llVert  u_{K^{c}}\rrVert _{1}\leq\alpha\llVert
u_{K}\rrVert _{1} \}$.
We may want to control the size of selected models with the upper
sparse eigenvalue, defined as
\[
\kappa^*(m,\widetilde{\mathbf X}) = \max_{\|u\|=1,\|u\|_0\le m}\|
\widetilde{\mathbf{X}} u\|^2/n. 
\]
%

We impose the following conditions on the {target} coefficient vector and
the design:
%
\begin{equation}
\label{cond-4a} \mathbb{P} \{\operatorname{Cond}_1 \}\ge1-{\widetilde
\epsilon}_1, \qquad\operatorname{Cond}_1 = \biggl\{|K| + \sum
_{k\notin K}\frac{|{\bolds
{\gamma
}}^{\mathrm{target}}_k/\sigma^{\mathrm{ora}}|}{
\sqrt{(2/n)\log{\widetilde p}}} \le s_1 \biggr
\}
\end{equation}
%
for a certain index set $K$, and
%
\begin{equation}
\label{cond-5a} \mathbb{P} \{\operatorname{Cond}_2 \}\ge1-{\widetilde
\epsilon}_1, \qquad\operatorname{Cond}_2 = \Bigl\{ \max
_{|J\setminus K| \le s_2} \phi ^{-2}_{\mathrm{comp}} ( \alpha,J,
\widetilde{\mathbf X} ) \le C_2 \Bigr\}.
\end{equation}
%
For small penalty levels and the LSE after model selection, we also need
%
\begin{equation}
\label{cond-6a} \mathbb{P} \{\operatorname{Cond}_3 \}\ge1-{
\widetilde\epsilon}_1, \qquad\operatorname{Cond}_3 = \bigl\{
\kappa^*(s_3,\widetilde{\mathbf X}) \le C_3 \bigr\}.
\end{equation}
Finally, for ${\bolds{\gamma}}^{\mathrm{target}}\neq\bolds{\gamma}$, we need
the condition
%
\begin{eqnarray}\label{cond-7a}
\mathbb{P} \{\operatorname{Cond}_4 \}&\ge&1-{\widetilde
\epsilon}_1,
\nonumber
\\[-8pt]
\\[-8pt]
\nonumber
 \operatorname {Cond}_4 &=& \bigl\{
C_4\bigl\|\widetilde{\mathbf X}\bigl({\bolds{\gamma }}^{\mathrm{target}}-
\bolds {\gamma}\bigr)\bigr\| \le\sigma^{\mathrm{ora}}\sqrt{\log({\widetilde p}/{
\widetilde\epsilon }_1)} \bigr\}. 
\end{eqnarray}
In (\ref{cond-4a}), (\ref{cond-5a}), (\ref{cond-6a}) and (\ref{cond-7a}),
$s_j$ are allowed to change with $\{n,{\widetilde p}\}$, while $\alpha$
and $C_j$ are fixed constants.
These conditions also make sense for deterministic designs with
${\widetilde\epsilon}_1=0$ for
deterministic conditions.

Let $k$ and $\varepsilon$ be positive real numbers and $\lambda_0$ be a
penalty level satisfying
%
\begin{equation}
\label{cond-8a} \lam_0\ge(1+\varepsilon) L_{n-3/2}(k/{
\widetilde p}), 
\end{equation}
where $L_n(t) = n^{-1/2}\Phi^{-1}(1-t)$ is the $\mathcal{N}(0,1/n)$
negative quantile function. Let
%
\begin{equation}
\varepsilon_1 \ge\frac{e^{1/(4n-6)^2}4k/{{s_2}}}{L_1^4(k/{\widetilde
p})+2L_1^2(k/{\widetilde p})} + \biggl(\frac{L_1({\widetilde\epsilon}_1/{\widetilde
p})}{L_1(k/{\widetilde p})}+
\frac{e^{1/(4n-6)^2}/\sqrt{2\pi}}{L_1(k/{\widetilde p})} \biggr)\sqrt {\frac
{C_3}{{{s_2}}}}. \label{cond-9a}
\end{equation}
We note that $L_n(t)=n^{-1/2}L_1(t)\le\sqrt{(2/n)\log(1/t)}$ for
$t\le
1/2$, so that the right-hand
side of (\ref{cond-9a}) is of the order $k/\{s_2(\log p)^2\}+\sqrt
{\delta/s_2}$. Thus condition
(\ref{cond-9a}) is easily satisfied even when $\varepsilon_1$ is a
small positive number
and $k$ is a moderately large number. Moreover, $\lam$ depends on
$\delta$ only through
$\sqrt{\delta/s_2}$ in (\ref{cond-9a}).

%
\begin{theorem}\label{th-lasso} Let $\{\hat{\bolds \gamma},\hat
{\sigma
}\}$ be as in
(\ref{standard Scaled Lasso}) with data in (\ref{LM-a}) and a penalty
level in (\ref{cond-8a}).
Let ${\widetilde\epsilon}_1<1$ and $\lam^*=L_{n-3/2}({\widetilde
\epsilon}_1/{\widetilde p})$.
Suppose $\lam^*\le1$ and $\delta s(\log{\widetilde p})/\break n \le c_0$.

%
\begin{longlist}[(iii)]
\item[(i)] Let ${\bolds{\gamma}}^{\mathrm{target}} = \bolds{\gamma}$,
$s\ge
s_1$, $s_2=0$,
$k\le{\widetilde\epsilon}_1$ in (\ref{cond-8a}), $\alpha
=1+2/\varepsilon$
and $4{\widetilde\epsilon}_1\le{\widetilde p}^{1-\delta}{\widetilde
\epsilon}_0$.
Then there exists a constant $C_0$ depending on $\{\alpha,C_2\}$ only
such that
when $C_0c_0\le1/2$,
(\ref{cond-4a}) and (\ref{cond-5a}) imply (\ref{cond-1a}), (\ref
{cond-2a}) and (\ref{cond-3a}).
%
\item[(ii)] Let ${\bolds{\gamma}}^{\mathrm{target}} = \bolds{\gamma}$,
$s\ge
s_1+{{s_2}}$, $1\le s_2\le s_3$,
$k\ge1$ and $\varepsilon_1<\varepsilon$ in (\ref{cond-8a}) and
(\ref{cond-9a}),
$\alpha\ge\sqrt{2}(\varepsilon-\varepsilon_1)_+^{-1} \{
1+\varepsilon+L_{1}({\widetilde\epsilon}_1/{\widetilde
p})/L_{1}(k/{\widetilde p}) \}$
and $(5+e^{1/(4n-6)^2}){\widetilde\epsilon}_1\le{\widetilde
p}^{1-\delta}{\widetilde\epsilon}_0$.
Then there exists a constant $C_0$ depending on $\{\alpha,\varepsilon
,\varepsilon_1,C_2\}$ only such that
when $C_0c_0\le1/2$, (\ref{cond-4a}), (\ref{cond-5a}) and (\ref
{cond-6a}) imply (\ref{cond-1a}),
(\ref{cond-2a}) and (\ref{cond-3a}).
\item[(iii)] Let $s\ge s_1+{{s_2}}$, $1\le s_2\le s_3$,
$k\ge1$ and $\varepsilon_1<\varepsilon$ in (\ref{cond-8a}) and
(\ref{cond-9a}),
$\varepsilon_1<\varepsilon_2<\varepsilon$,
$\alpha\ge2(\varepsilon-\varepsilon_2)_+^{-1} \{
1+\varepsilon+L_{1}({\widetilde\epsilon}_1/{\widetilde
p})/L_{1}(k/{\widetilde p}) \}$,
$(6+e^{1/(4n-6)^2}){\widetilde\epsilon}_1\le{\widetilde p}^{1-\delta
}{\widetilde\epsilon}_0$
and $C_4 \ge\sqrt{(4/L_1^2(k/{\widetilde p}))\log({\widetilde
p}/{\widetilde\epsilon}_1)}/
\min(\sqrt{2}-1,\varepsilon_2-\varepsilon_1)$.
Then there exists a constant $C_0$ depending on $\{\alpha,\varepsilon
,\varepsilon_2,C_2\}$ only such that
when $C_0c_0\le1/2$, (\ref{cond-4a}), (\ref{cond-5a}), (\ref{cond-6a})
and (\ref{cond-7a}) imply (\ref{cond-1a}), (\ref{cond-2a}) and (\ref
{cond-3a}).
\end{longlist}
\end{theorem}

In Theorem~\ref{th-lasso}, $s_1$ in (\ref{cond-4a}) represents the
complexity or the size of
the coefficient vector, and $s_2$ represents the number of false
positives we are
willing to accept with the penalty level in (\ref{cond-8a}). Thus $s$
is an upper bound
for the total number of estimated coefficients, true or false.
We summarize parallel results for the LSE after model selection as follows.

\begin{theorem}
\label{Error bound of LSE}
Let $\{\hat{\bolds \gamma},\hat{\sigma}\}$ be as in (\ref{standard
Scaled Lasso}) and
$ \{\hat{\bolds \gamma}^{\mathrm{lse}},\hat{\sigma}^{\mathrm{lse}} \}$
as in (\ref{Scaled Lasso LSEa}).
\begin{longlist}[(ii)]
\item[(i)] The following bounds always hold:
%
\begin{equation}
{\hat\sigma}^2 - \bigl({\hat\sigma}^{\mathrm{lse}}
\bigr)^2 =\bigl \|\widetilde{\mathbf X} \bigl(\hat{\bolds
\gamma}^{\mathrm{lse}} - \hat{\bolds\gamma} \bigr) \bigr\|^2/n \le
\frac{({\hat\sigma}\lam_0)^2|{\hat S}|}{
\phi_{\mathrm{comp}}^{2} (0,{\hat S},\widetilde{\mathbf  X} )} \label{LSE 1}
\end{equation}
with ${\hat S}=\operatorname{ supp}(\hat{\bolds \gamma})$ and
%
\begin{equation}
\bigl\|\hat{\bolds \gamma}^{\mathrm{lse}} - \hat{\bolds \gamma}\bigr\|_1 \le
\frac{{\hat\sigma}\lam_0|{\hat S}|}{\phi_{\mathrm
{comp}}^{2}
(0,{\hat S},\widetilde{\mathbf  X} )}. \label{LSE 2}
\end{equation}
\item[(ii)] Let $\lam_0$ be a penalty level satisfying (\ref
{cond-8a}) and
$\varepsilon_1<\varepsilon_2<\varepsilon_3<\varepsilon$.
Suppose the conditions of Theorem~\ref{th-lasso} hold and that the
constant factor $C_0$
in Theorem~\ref{th-lasso} satisfies
\begin{eqnarray*}
C_0s\delta(\log{\widetilde p})/n\le\frac{\varepsilon-\varepsilon
_3}{1+\varepsilon},\qquad
\frac{C_0s\delta(\log{\widetilde p})}{(\varepsilon_3-\varepsilon
_2)^2L_1^2(k/{\widetilde p})} \le\frac{s_3}{C_3}.
\end{eqnarray*}
Then, for the parameters defined in the respective parts of Theorem~\ref
{th-lasso},
%
\begin{equation}
\mathbb{P} \bigl\{ \llvert {\hat S}\rrvert < s_3+s \bigr\} \ge1- {
\widetilde p}^{1-\delta}{\widetilde\epsilon}_0. \label{dim bound}
\end{equation}
If in addition, condition (\ref{cond-5a}) is strengthened to
%
\begin{equation}
\mathbb{P} \Bigl\{ \Bigl(\max_{|J\setminus K| \le s_2} \phi
^{-2}_{\mathrm{comp}} ( \alpha,J,\widetilde{\mathbf X} ) \Bigr)\vee
\Bigl( \max_{|J| \le s_3+s} \phi^{-2}_{\mathrm{comp}} (0,J,
\widetilde{\mathbf X} ) \Bigr) > C_2 \Bigr\}\le{\widetilde
\epsilon}_1, \label{cond-5b}
\end{equation}
then the conclusions of Theorem~\ref{th-lasso} hold with $\{\hat{\bolds\gamma},\hat{\sigma}\}$ replaced by
$ \{\hat{\bolds \gamma}^{\mathrm{lse}},\hat{\sigma}^{\mathrm{lse}} \}$.
\end{longlist}
\end{theorem}

We collect some probability bounds for the regularity conditions in the
following proposition.
Consider deterministic coefficient vectors $\bolds{\beta}^{\mathrm{target}}$
satisfying
%
\begin{equation}
\label{cond-4b} |K| + \sum_{j\notin K} \frac{C_1|\beta_j^{\mathrm{target}}|}{\sqrt
{(2/n)\log
{\widetilde p}}}
\le s_1.
\end{equation}

\begin{proposition}\label{prop-probab}
Let $\mathbf{X}$ be a $n\times p$ matrix with i.i.d. $\mathcal
{N}(0,\bolds{\Sigma})$ rows,
$A^c\subset\{1,\ldots,p\}$ with $|A^c|={\widetilde p}$,
$\Dbar= \operatorname{ diag}(\mathbf{X}^T\mathbf{X}/n)$,
$\widetilde{\mathbf  X}=\mathbf{X}_{A^c}\Dbar_{A^c}^{-1/2}$,
$\bolds{\gamma} = \Dbar_{A^c}^{1/2}\bolds{\beta}_{A^c}$ and
$\bolds{\gamma}^{\mathrm{target}} = \Dbar_{A^c}^{1/2}\bolds{\beta
}_{A^c}^{\mathrm{target}}$.
Suppose $1/M\le\lam_{\min}(\bolds{\Sigma})\le\lam_{\max
}(\bolds
{\Sigma})\le M$ with a fixed $M$.
Let $\lam_1 = \sqrt{(2/n)\log({\widetilde p}/{\widetilde\epsilon}_1)}$.
Then, for a certain constant $C_*$ depending on $M$ only,
%
\begin{eqnarray}
\label{probab-1a} & &(\ref{cond-4b})\Rightarrow(\ref{cond-4a}) \qquad\mbox{when }
C_1\ge \sqrt{M(1+ \lam_1)},
\\
\label{probab-2a} && C_2\ge C_*\bigl\{1+ \max \bigl\{|K|+s_2,s+s_3
\bigr\} \lam_1^2\bigr\} \Rightarrow (\ref{cond-5b})
\Rightarrow(\ref{cond-5a}),
\\
\label{probab-3a} && C_3\ge C_*\bigl\{1+ s_3
\lam_1^2\bigr\}\Rightarrow (\ref
{cond-6a}),
\end{eqnarray}
and for $\lam_2 = \sqrt{(2/n)\log(1/{\widetilde\epsilon}_1)}$ and
any coefficient vectors $\bolds{\beta}^{\mathrm{target}}$ and $\bolds{\beta}$,
%
\begin{equation}
\label{probab-4a} C_4C_*(1+\lam_2)\bigl\|\bolds{
\beta}^{\mathrm{target}}-\bolds{\beta}\bigr\| \le \lam_1 \Rightarrow
(\ref{cond-7a}).
\end{equation}
Moreover, when $\lam_2\le1/2$, $\sigma^{\mathrm{ora}}$ can be replaced by
$\sqrt{\mathbb{E}(\sigma^{\mathrm{ora}})^2}$
or $C_*$ in (\ref{cond-1a}) and~(\ref{cond-2a}). \
\end{proposition}

\section{Discussion}
\label{Discu0}

\subsection{Alternative choice of penalty level for finite sample performance}
\label{Discu}

In Theorem~\ref{Main} and nearly all consequent results in
Theorems \ref{Main1}--\ref{Main 2} and \ref{support recovery}--\ref%
{support recovery Latent}, we have picked the penalty level $\lambda
= (
1+\varepsilon )\sqrt{(2\delta/n) \log p}$ for $\delta\geq
1 $ ($\delta\geq3$ for support recovery) and $\varepsilon>0$. This
choice of $\lambda$ can be too conservative and may cause some finite sample
estimation bias. However, in view of Theorem~\ref{th-lasso}(ii) and (iii),
the results in these theorems in Sections~\ref{Methodology} and \ref{app} still hold for
penalty levels no smaller
than $\lam=  (1+\varepsilon )L_{n}(k/p)
\approx (1+\varepsilon )\sqrt{(2/n)\log(p/k)}$,
which weakly depends on $\delta$ through (\ref{cond-9a})
and the requirement of $\varepsilon>\varepsilon_1$.

Condition (\ref{cond-9a}), with $\varepsilon<\varepsilon_1$,
${\widetilde\epsilon}_1={\widetilde p}^{1-\delta}$
and ${\widetilde p}=p-2$ for the estimation of precision matrix,
is the key for the choice of the smaller penalty level $\lam=
(1+\varepsilon )L_{n}(k/p)$.
%
It provides theoretical justifications for the choice of $k\in[1, n]$
or even up to $k \asymp n\log p$
for the theory to work.
Let $s_{\max} = c_{0}n/\log p$ with a sufficiently small constant
$c_{0}>0$, which can be viewed as the largest possible $s\ge s_1+s_2$
in our theory.
Suppose $n \le p^{t_0}$ for some fixed $t_0<1$ and the bound $C_3$ for the
upper sparse eigenvalue can be treated as fixed in (\ref{cond-6a}) for
$s_2\le s_{\max}$.
For $\lam=  (1+\varepsilon )\sqrt{(2/n)\log(p/k)}$ with
$k\le n\log p$ and $s_2\le s_{\max}$,
condition (\ref{cond-9a}) can be written as
\[
\varepsilon> \varepsilon_1 \ge\frac{(1+o(1))(k/n)s_{\max
}/(c_0s_2)}{(1-t_0)(1+(1-t_0)\log p)} +\bigl(\sqrt{
\delta}+o(1)\bigr)\sqrt{C_3/s_2},
\]
which holds for sufficiently small $k s_{\max}/(ns_2\log p)$.
This allows $k\asymp n\log p$ for $s_2=s_{\max}$.
For the asymptotic normality, we need $s_2=o(\sqrt{n}/\log p)$, so that
$k=o(n^{1/2}\log p)$ is sufficient.

\subsection{Statistical inference under unbounded condition number}

The main results in this paper assume that the spectrum of the precision
matrix $\Omega$ is bounded below and above by a universal constant $M$ in
(\ref{sparseparaspace}). Then the dependency of the key result
(\ref{result1}) on $M$ is hidden in the constant $C_{1}$ in front
of the rate $s\frac{\log p}{n}$ for bounding $\llVert \hat{\theta
}%
_{A,A}-\theta_{A,A}^{\mathrm{ora}}\rrVert _{\infty}$ in Theorem~\ref
{Main}. The
inference result in (\ref{result 4a}) follows as long as this bound
$C_{1}s\frac{\log p}{n}$ is
dominated by the parametric square-root rate of $\theta_{A,A}^{\mathrm{ora}}$,
or equivalently $s=o (\sqrt{n}/\log p ) $. In
fact, following the proof of Theorem~\ref{Main}, the requirement can be
somewhat weakened to $\lambda_{\max}(\Omega)\leq M$ and $\max_{j}\sigma_{jj}\leq M$.

It would be interesting to consider a slightly more general case, where we
assume $\max_{j}\sigma_{jj}\leq C$ and $\lambda_{\max}(\Omega)\leq
M_{n,p}$ with absolute constant $C$ and possibly a large constant $%
M_{p,n}\rightarrow\infty$ as ($n,p$)$\rightarrow\infty$. In this
setting, the condition number of $\Omega$ may not be bounded. Suppose we
would like to make inference for $\omega_{12}$ and assume $\max
\{
\omega_{11},\omega_{22} \} \leq C$ to make its inverse Fisher
information bounded. We are able to show that (\ref{result 4a}) holds
as long
as $s=o ( \sqrt{n}/ ( M_{n,p}\log p )  ) $ under this
setting. In fact, the regression model (\ref{CondiDis2}) is still
valid with
bounded noise level $\theta_{mm}\leq\sigma_{mm}\leq C$ for $m\in
A= \{ 1,2 \} $. However, the compatibility condition (\ref{cond-5a})
may not hold with absolute constant because the smallest
eigenvalue of the population Gram matrix $\lambda_{\min}(\Sigma
_{A^{c}A^{c}})$ is possibly as small as $M_{n,p}^{-1}$. Taking this possible
compatibility constant $M_{n,p}^{-1}$ into account, we can obtain $%
\llVert \hat{\theta}_{A,A}-\theta_{A,A}^{\mathrm{ora}}\rrVert
_{\infty
}=O_{p}(sM_{n,p}\frac{\log p}{n})$ while the sufficient statistics
$\theta
_{A,A}^{\mathrm{ora}}$ still has square-root rate. As a consequence, the inference
result in (\ref{result 4a}) holds as long as $s=o ( \sqrt{n}%
/ ( M_{n,p}\log p )  ) $. We would like to point out
that to
guarantee compatibility condition (\ref{cond-5a}) indeed holds at
the level $C_2\asymp M_{n,p}$, an extra condition $\sqrt{s(\log p)/n}%
=o(M_{n,p}^{-1})$ is required; see Corollary~1 in \citet%
{raskutti2010restricted}. However, when $s=o ( \sqrt{n}/ (
M_{n,p}\log p )  ) $, this condition is automatically satisfied.
The argument above can be made rigorous.

\subsection{Related works}
\label{Discu1}

Our methodology in this paper is related to \citet{zhang2014confidence} who
proposed a LDPE approach for making inference in a high-dimensional
linear model.
Since $\hat{\bolds {\epsilon}}_{A}$ can be viewed as an approximate
projection of $\mathbf{X}_{A}$ to the direction of $\bolds{\epsilon
}_{A}$ in (\ref{estimator}),
the estimator in (\ref{estimator}) can be viewed as an LDPE as \citet
{zhang2014confidence}
discussed in the regression context.
See also \citet{GeerBR13} and \citet{JavanmardM13}.
When appropriately applying their approach to our setting,
their result is asymptotically equivalent to ours and also obtains the
asymptotic normality. In this section, we briefly discuss their
approach in
the large graphical model setting.

Consider $A=\{1,2\}$.
While our method regresses two nodes $\mathbf{X}_{A}$ against all other
nodes $\mathbf{X}_{A^{c}}$ and focuses on the estimation of
the two by two dimensional covariance matrix $\Omega_{A,A}^{-1}$ of the
noise, their approach consists of the following two steps. First, one
node $%
\mathbf{X}_{1}$ is regressed against all other nodes $\mathbf{X}_{1^{c}}$
using scaled lasso with coefficient $\hat{\beta}^{\mathrm{(init)}}$. As equation
(%
\ref{CondiDis2}) suggests, the noise level is $\omega_{11}^{-1}$, and the
coefficient for the column $\mathbf{X}_{2}$ is $\beta_{2}=-\omega
_{12}\omega_{11}^{-1}$. Then in the second step, to correct the bias
of the
initial estimator $\hat{\beta}_{2}^{\mathrm{(init)}}$ obtained in the first
step for
the coefficient vector $\beta_{2}$, a score vector ${\mathbf z}$ is
picked and applied to
obtain the final estimator of $\beta_{2}$ as follows:
\[
\hat{\beta}_{2}=\hat{\beta}_{2}^{\mathrm{(init)}}+{\mathbf
z}^{T}\bigl(\mathbf {X}_{1}-\mathbf{X}%
_{1^{c}}
\hat{\beta}_{2}^{\mathrm{(init)}}\bigr)/{\mathbf z}^{T}
\mathbf{X}_{2},
\]
where ${\mathbf z}$ is the residue after regressing $\mathbf{X}_{2}$
against all
remaining columns in step one $\mathbf{X}_{A^{c}}$ using scaled lasso again.
To obtain the final estimator of $\omega_{12}$, the estimator $\hat
{\beta}%
_{2}$ of $-\omega_{12}\omega_{11}^{-1}$ should be scaled by an accurate
estimator of $\omega_{11}^{-1}$, which uses the variance component of
the scaled lasso
estimator in the first step. It seems that two approaches are quite different.
However, both approaches do the same thing: they try to estimate the partial
correlation of node $Z_{1}$ and $Z_{2}$ and hence are asymptotically
equivalent. Compared with their approach, our method enjoys simper form and
clearer interpretation. It is worthwhile to point out that the main
contribution of this paper is understanding the fundamental limit of the
Gaussian graphical model in making statistical inference, which is not
covered by other works.

\subsection{Unknown mean \texorpdfstring{$\mu$}{mu}}

In the \hyperref[sec1]{Introduction}, we assume $Z\sim\mathcal{N} ( \mu,\Sigma
 ) $
and $\mu=0$ without loss of generality. This can be seen as follows.
Suppose we observe an $n\times p$ data matrix $\mathbf{X}$ with i.i.d.
rows from
$\mathcal{N} ( \mu,\Sigma )$.
Let $u^{(i)}$, $i=1,\ldots,n$, be $n$-dimensional orthonormal row
vectors with $u^{(n)}=(1,\ldots,1)/\sqrt{n}$.
Then $u^{(i)}\mathbf{X}$ are i.i.d. $p$-dimensional row vectors from
$\mathcal{N} (0 ,\Sigma )$.
Thus we can simply apply our methods and theory to the sample $\{
u^{(i)}\mathbf{X}, i =1,\ldots,n-1\}$.

\section{Numerical studies}
\label{Numer}

In this section, we present some numerical results for both asymptotic
distribution and support recovery. We generate the data from $p\times p$
precision matrices with three blocks. Two cases are considered: $p=200,800$.
The ratio of block sizes is $2:1:1$; that is, for a $200\times200$ matrix,
the block sizes are $100\times100$, $50\times50$ and $50\times50$,
respectively.
The diagonal entries are $\alpha_{1},\alpha_{2},\alpha_{3}$ in three
blocks, respectively,
where $(\alpha_{1},\alpha_{2},\alpha_{3})=(1,2,4)$.
When the entry is in the $k$th block, $\omega_{j-1,j}=\omega
_{j,j-1}=0.5\alpha_{k}$, and $\omega
_{j-2,j}=\omega_{j,j-2}=0.4\alpha_{k}$, $k=1,2,3$. The asymptotic variance
for estimating each entry can be very different. Thus a simple
procedure with a single
threshold level for all entries is not likely to perform well.

We first estimate the entries in the precision matrix and partial
correlations as discussed in Remark~\ref{Partial Correlation}, and
consider the distributions of these estimators. We generate a random
sample of
size $n=400$ from a multivariate Gaussian distribution $\mathcal
{N}(0,\Sigma)$
with $\Sigma=\Omega^{-1}$. For the proposed estimators defined through
(\ref{estimator}) and (\ref{estimator1}) with the scaled lasso (\ref
{Scaled Lasso})
or the LSE after model selection (\ref{Scaled Lasso LSE}), we pick
$\lam= n^{-1/2}L_n(1/p) \approx\sqrt{(2/n)\log p}$; that is, $k=1$ in
(\ref{cond-8a})
with small adjustment in $n$ and $p$ ignored. This is justified by our
theoretical
results as discussed in Section~\ref{Discu}.


\begin{table}
\caption{Mean and standard error of GLasso, CLIME and proposed estimators}
\label{table-infer}
\begin{tabular*}{\textwidth}{@{\extracolsep{\fill}}lccccc@{}}
\hline
$\mathbf{p}$& & $\bolds{\omega_{1,2}=0.5}$ & $\bolds{\omega_{1,3}=0.4}$ &
$\bolds{\omega_{1,4}=0}$ &
$\bolds{\omega_{1,10}=0}$  \\
\hline
200
& GLasso & 0.368 $\pm$ 0.039& 0.282 $\pm$ 0.038& $-$0.056 $\pm$ 0.03\phantom{0}&
$-$0.001 $\pm$ 0.01\phantom{0}\\
& CLIME & 0.776 $\pm$ 0.479& 0.789 $\pm$ 0.556&\phantom{0} 0.482 $\pm$ 1.181&
\phantom{0.}0.002 $\pm$ 0.017\\
& $\hat{\omega}_{i,j}$ & 0.459 $\pm$ 0.05\phantom{0}& 0.372 $\pm$ 0.052&
$-$0.049 $\pm$ 0.041& $-$0.003 $\pm$ 0.044\\[2pt]
& $\hat{\omega}_{i,j}^{\mathrm{LSE}}$ & 0.503 $\pm$ 0.059& 0.401
$\pm$ 0.061& $-$0.006 $\pm$ 0.049& $-$0.002 $\pm$ 0.052\\[3pt]
800
& GLasso & 0.801 $\pm$ 0.039& 0.258 $\pm$ 0.031& \phantom{00.}0.19 $\pm$ 0.014&
$-$0.063 $\pm$ 0.028\\
& CLIME & 1.006 $\pm$ 0.255& 0.046 $\pm$ 0.140 & \phantom{0.}0.022 $\pm$ 0.071&
\phantom{0.}0.018 $\pm$ 0.099\\
& $\hat{\omega}_{i,j}$ & 0.436 $\pm$ 0.049& 0.361 $\pm$ 0.047&
$-$0.057 $\pm$ 0.044& \phantom{0.}0.001 $\pm$ 0.044\\[2pt]
& $\hat{\omega}_{i,j}^{\mathrm{LSE}}$ & 0.491 $\pm$ 0.059& 0.396
$\pm$ 0.058& \phantom{00000}0 $\pm$ 0.052& $-$0.003 $\pm$ 0.05\phantom{0}\\
\hline
$\mathbf{p}$ & & $\mathbf{r_{1,2}=-0.5}$ & $\mathbf{r_{1,3}=-0.4}$ &
$\mathbf{r_{1,4}=0}$ & $\mathbf{r_{1,10}=0}$ \\
\hline
200 & $\hat{r}_{i,j}$ & $-$0.477 $\pm$ 0.037& $-$0.391 $\pm$ 0.043&
0.051 $\pm$ 0.043& \phantom{0.}0.003 $\pm$ 0.046\\[2pt]
& $\hat{r}_{i,j}^{\mathrm{LSE}}$ & $-$0.485 $\pm$ 0.04\phantom{0}& $-$0.386 $\pm$
0.046& 0.006 $\pm$ 0.047& \phantom{0.}0.002 $\pm$ 0.049\\[3pt]
800 & $\hat{r}_{i,j}$ & $-$0.468 $\pm$ 0.039& $-$0.392 $\pm$ 0.041&
\phantom{0}0.06 $\pm$ 0.045& $-$0.001 $\pm$ 0.048\\[2pt]
& $\hat{r}_{i,j}^{\mathrm{LSE}}$ & $-$0.475 $\pm$ 0.041& $-$0.382
$\pm$
0.044& \phantom{000.}0 $\pm$ 0.049& \phantom{0.}0.002 $\pm$ 0.048\\
\hline
\end{tabular*}
\end{table}

\begin{figure}

\includegraphics{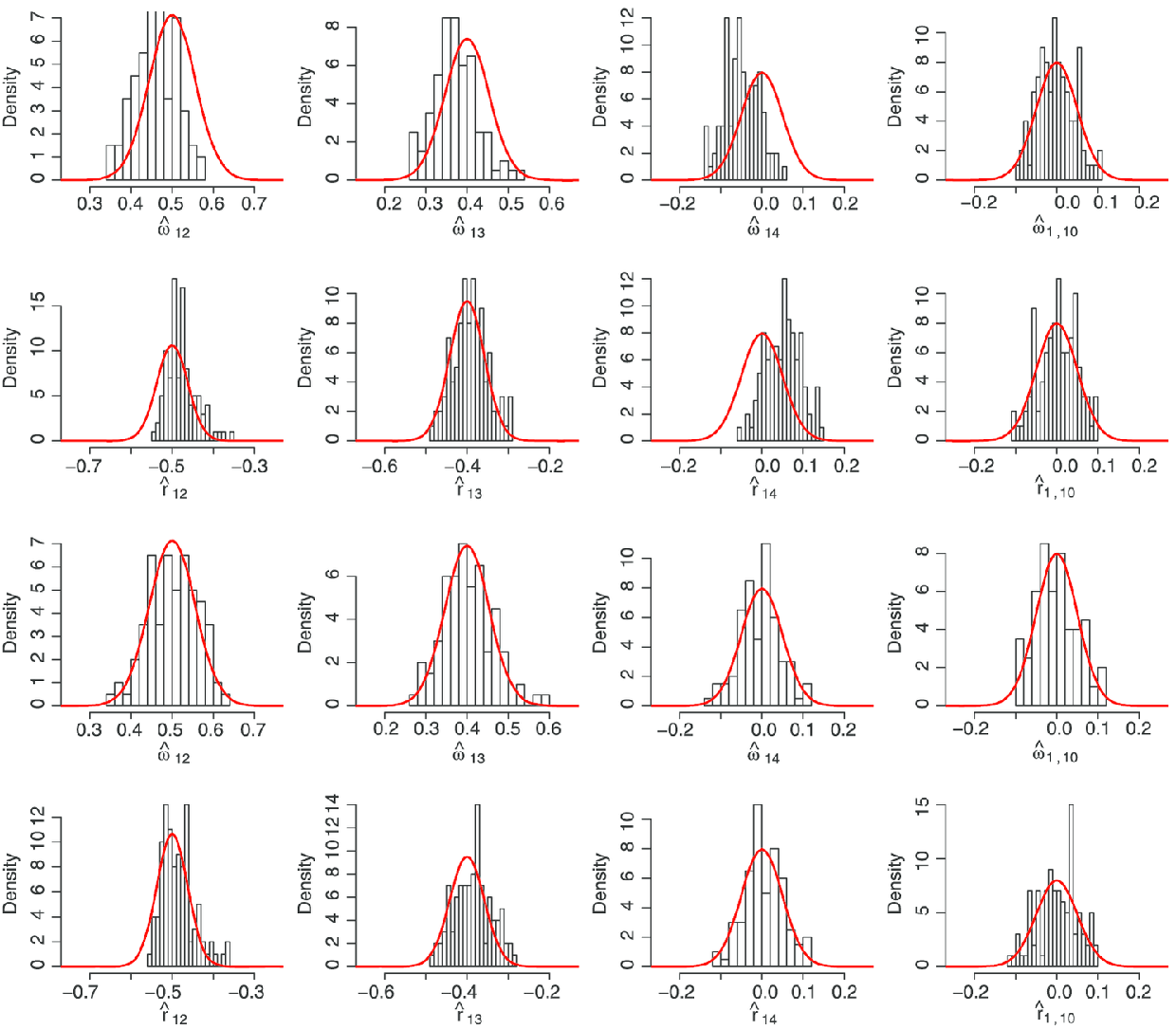}

\caption{Histograms of estimated entries for $p=200$.
The first row: scaled lasso for entries $\protect\omega
_{1,2} $ and $\protect\omega_{1,3}$ in the precision matrix;
the second row: scaled lasso for entries
$\protect\omega_{1,4}$ and $\protect\omega_{1,10}$;
the third and fourth rows: LSE after scaled lasso selection.}\label
{fig-hist-p200}
\end{figure}

\begin{figure}

\includegraphics{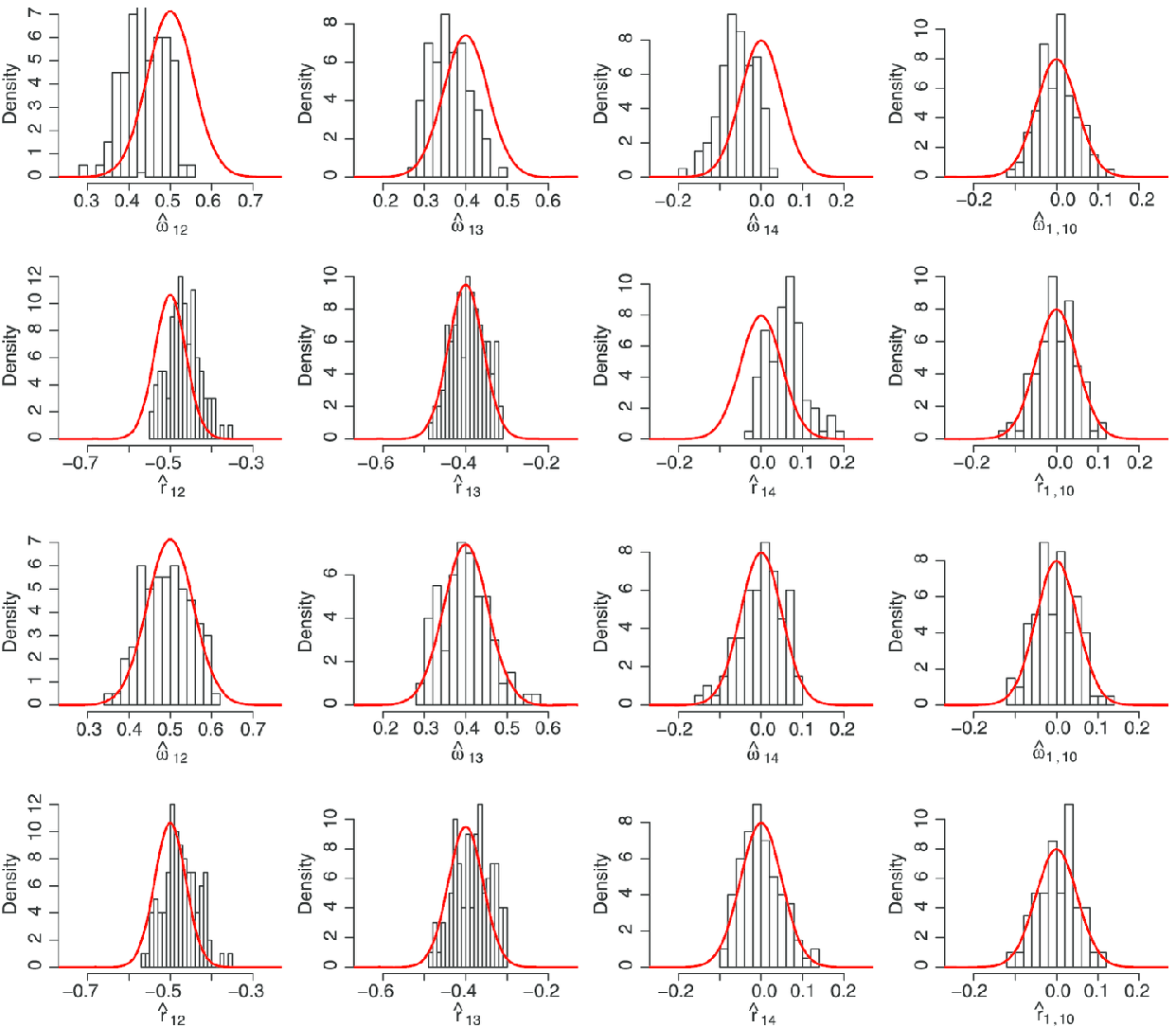}

\caption{Histograms of estimated entries for $p=800$.
The first row: scaled lasso for entries $\protect\omega
_{1,2} $ and $\protect\omega_{1,3}$ in the precision matrix;
the second row: scaled lasso for entries
$\protect\omega_{1,4}$ and $\protect\omega_{1,10}$;
the third and fourth rows: LSE after scaled lasso selection.}
\label{fig-hist-p800}\vspace*{-5pt}
\end{figure}

Table~\ref{table-infer} reports the mean and standard error of our
estimators for four entries in the precision matrices and the corresponding
correlations. In addition, we report the point estimates by the GLasso
[\citet{FHT08}] and CLIME [\citet{CLL11}] for comparison.
For $p=800$, the results for the GLasso are based on 10 replications,
while all
other entries in the table are based on 100 replications.
The GLasso is computed by the R~package ``glasso'' with penalized
diagonal (default option),
while the CLIME estimators are computed by the R package ``fastclime''
[\citet{JMLR:v15:pang14a}].
As the GLasso and CLIME are designed for estimating precision matrices
as high-dimensional objects,
it is not surprising that the proposed estimator outperforms them in
estimation accuracy
for individual entries.
Figures~\ref{fig-hist-p200} and \ref{fig-hist-p800} show the histograms
of the proposed
estimates with the theoretical Gaussian density in Theorem~\ref{Main1}
super-imposed.
They demonstrated that the histograms match pretty well to the
asymptotic distribution,
especially for the LSE after model selection.
The asymptotic normality leads to the following $(1-\alpha)$ confidence
intervals
for $\omega_{ij}$ and $r_{ij}$:
\begin{eqnarray*}
&&\Bigl(\hat{\omega}_{ij} - z_{\alpha/2}\sqrt{\bigl(
\hat{\omega }_{i,i}\hat{\omega}_{j,j}+\hat{
\omega}_{i,j}^2\bigr)/n},  \omega_{ij} +
z_{\alpha/2}\sqrt{\bigl(\hat{\omega}_{i,i}\hat {
\omega }_{j,j}+\hat{\omega}_{i,j}^2\bigr)/n}
\Bigr),
\\
&&\bigl(\hat{r}_{i,j} - z_{\alpha/2}\bigl(1-
\hat{r}_{i,j}^2\bigr)/\sqrt {n},  \hat{r}_{i,j}
+ z_{\alpha/2}\bigl(1-\hat{r}_{i,j}^2\bigr)/\sqrt
{n} \bigr),
\end{eqnarray*}
where $z_{\alpha/2}$ is the $z$-score such that $P(\mathcal
{N}(0,1)>z_{\alpha/2})=\alpha/2$.
Table~\ref{table-interval} reports the empirical coverage probabilities
for 95\% confidence intervals, which matches well to the assigned
confidence level.

Support recovery of a precision matrix is of great interest. We compare our
selection results with the GLasso and CLIME. In addition to the
training sample, we
generate an independent sample of size 400 from the same distribution for
validating the tuning parameter for the GLasso and CLIME. These
estimators are
computed based on the entire training sample with a range of penalty
levels and
a proper penalty level is chosen by minimizing the negative likelihood
$\{\operatorname{trace}(%
\overline{\Sigma}\hat{\Omega})-\log\det(\hat{\Omega})\}$
on the
validation sample, where $\overline{\Sigma}$ is the sample covariance matrix.
The proposed ANT estimators are computed based on the training sample only
with $\xi_{0}=2$ in the thresholding step as in (\ref{thresholding level}).
Tables~\ref{table-supp} and \ref{table-supp-800} present the average selection
performances as measured in the true positive, false positive and the
corresponding rates.
In addition to the overall performance, the summary statistics are
reported for each block.
The results demonstrate the selection consistency property of both ANT methods
and substantial false positive for the GLasso and CLIME.
It should be pointed out that the ANT takes the advantage of an
additional thresholding
step, while the GLasso and CLIME do not.
A possible explanation of the false positive for the GLasso is a
tendency for the likelihood criterion with
the validation sample to pick a small penalty level.
However, such an explanation seems not to hold for the CLIME, which
demonstrated much lower
false positive than the GLasso, as the true positive rate of the CLIME
is consistently maintained at
about 95\% for $p=200$ and 85\% for $p=800$.

Moreover, we compare the ANT with the GLasso and CLIME in a range of
penalty levels.
Figure~\ref{fig-roc} plots the ROC curves for the GLasso and CLIME with
various penalty levels and
the ANT with various thresholding levels in the follow-up procedure.
It demonstrates that the CLIME outperforms the GLasso, but the two
methods perform
significantly more poorly than the ANT in the experiment.
In addition, the circle in the plot represents the performance of the
ANT with the selected threshold
level as in (\ref{thresholding level}). The triangle and diamond in the
plot represents the performance
of the GLasso and CLIME with the penalty level chosen by
cross-validation, respectively.
This again indicates that our method simultaneously achieves a very high
true positive rate and a very low false positive rate.

\section{Proof of Theorem \texorpdfstring{\protect\ref{sparseMinimaxlwSupport}}{5}}
\label{Thmproof5}

In this section we show that the upper bound given in Section~\ref{infer.sec}
is indeed rate optimal. We will only establish equation (\ref%
{sparserateSupport}). Equation~(\ref{sparserateSupport1}) is an immediate
consequence of equation (\ref{sparserateSupport}) and the lower bound
$\sqrt{%
\frac{\log p}{n}}$ for estimation of diagonal covariance matrices in
\citet%
{CZZH10}.

\begin{table}
\caption{Empirical coverage probabilities of the 95\% confidence intervals}
\label{table-interval}
\begin{tabular*}{\textwidth}{@{\extracolsep{\fill}}lccccc@{}}
\hline
$\mathbf{p}$ & $\mathbf{(i,j)}$ & \textbf{(1, 2)} & \textbf{(1, 3)} &
\textbf{(1, 4)} & \textbf{(1, 10)} \\
\hline
200 & $\hat{\omega}_{i,j}$  & 0.87& 0.89& 0.87& 0.98\\[2pt]
&$\hat{\omega}_{i,j}^{\mathrm{LSE}}$  & 0.96& 0.91& 0.94& 0.98\\
& ${\hat r}_{i,j}$  & 0.94& 0.94& 0.87& 0.98\\[2pt]
&${\hat r}_{i,j}^{\mathrm{LSE}}$  & 0.93& 0.94& 0.94& 0.97\\[3pt]
800 & $\hat{\omega}_{i,j}$  & 0.74& 0.88& 0.84& 0.95\\[2pt]
&$\hat{\omega}_{i,j}^{\mathrm{LSE}}$ & 0.93& 0.93& 0.96& 0.96\\
& ${\hat r}_{i,j}$  & 0.89& 0.98& 0.83& 0.95\\[2pt]
&${\hat r}_{i,j}^{\mathrm{LSE}}$  & 0.90& 0.94& 0.96& 0.96\\
\hline
\end{tabular*}
\end{table}

\begin{table}
\caption{The performance of support recovery ($p=200$, $100$ replications)}
\label{table-supp}
\begin{tabular*}{\textwidth}{@{\extracolsep{\fill}}lcd{3.2}d{1.4}d{4.2}d{1.4}@{}}
\hline
\textbf{Block} & \textbf{Method} & \multicolumn{1}{c}{\textbf{TP}} &
\multicolumn{1}{c}{\textbf{TPR}} & \multicolumn{1}{c}{\textbf{FP}} &
 \multicolumn{1}{c@{}}{\textbf{FPR}} \\
\hline
Overall
& GLasso & 391 & 1 & 5322.24 & 0.2728 \\
& CLIME & 372.61 & 0.953 & 588.34 & 0.0302 \\
& ANT & 391 & 1 & 0.04 & 0 \\
& ANT-LSE & 390.97 & 0.9999 & 0.01 & 0 \\[3pt]
Block 1 & GLasso & 197 & 1 & 1981.1 & 0.4168 \\
& CLIME & 188.47 & 0.9567 & 205.58 & 0.0433 \\
& ANT & 197 & 1 & 0 & 0 \\
& ANT-LSE & 196.98 & 0.9999 & 0 & 0 \\[3pt]
Block 2 & GLasso & 97 & 1 & 293.93 & 0.2606 \\
& CLIME & 92.29 & 0.9514 & 72.89 & 0.0646 \\
& ANT & 97 & 1 & 0 & 0 \\
& ANT-LSE & 96.99 & 0.9999 & 0 & 0 \\[3pt]
Block 3 & GLasso & 97 & 1 & 160.93 & 0.1427 \\
& CLIME & 91.85 & 0.9469 & 72.94 & 0.0647 \\
& ANT & 97 & 1 & 0 & 0 \\
& ANT-LSE & 97 & 1 & 0 & 0 \\
\hline
\end{tabular*}\vspace*{-7pt}
\end{table}
%
\begin{table}[b]\vspace*{-7pt}
\caption{The performance of support recovery ($p=800$, 10 replications)}
\label{table-supp-800}
\begin{tabular*}{\textwidth}{@{\extracolsep{\fill}}lcd{4.1}d{1.4}d{5.1}d{1.4}@{}}
\hline
\textbf{Block} & \textbf{Method} & \multicolumn{1}{c}{\textbf{TP}} & \multicolumn{1}{c}{\textbf{TPR}}
 & \multicolumn{1}{c}{\textbf{FP}} & \multicolumn{1}{c@{}}{\textbf{FPR}} \\
\hline
Overall
& GLasso & 1590.7 & 0.9998 & 44785.6 & 0.1408 \\
& CLIME & 1365.9 & 0.8585 & 134.6 & \multicolumn{1}{c}{4e--04} \\
& ANT & 1589 & 0.9987 & 0 & 0 \\
& ANT-LSE & 1586.2 & 0.997 & 0 & 0 \\[3pt]
Block 1
& GLasso & 797 & 1 & 19694.5 & 0.2493 \\
& CLIME & 687.5 & 0.8626 & 71.4 & \multicolumn{1}{c}{9e--04} \\
& ANT & 795.8 & 0.9985 & 0 & 0 \\
& ANT-LSE & 794.8 & 0.9972 & 0 & 0 \\[3pt]
Block 2
& GLasso & 397 & 1 & 2133.4 & 0.1094 \\
& CLIME & 339.4 & 0.8549 & 29.6 & 0.0015 \\
& ANT & 396.7 & 0.9992 & 0 & 0 \\
& ANT-LSE & 395.8 & 0.997 & 0 & 0 \\[3pt]
Block 3
& GLasso & 396.7 & 0.9992 & 664.7 & 0.0341 \\
& CLIME & 339 & 0.8539 & 32.6 & 0.0017 \\
& ANT & 396.5 & 0.9987 & 0 & 0 \\
& ANT-LSE & 395.6 & 0.9965 & 0 & 0 \\
\hline
\end{tabular*}\vspace*{-3pt}
\end{table}

The lower bound is established by Le Cam's method. To introduce Le Cam's
method we first introduce some notation. Consider a finite parameter
set $\mathcal{G}_{0}= \{ \Omega_{0},\Omega_{1},\ldots,\Omega
_{m_{\ast}} \} \subset\mathcal{G}_{0}(M,k_{n,p})$. Let
$\mathbb
{P}%
_{\Omega_{m}}$ denote the joint distribution of independent
observations $%
X^{(1)}$, $X^{(2)},\ldots, X^{(n)}$ with each $X^{(i)}\sim\mathcal
{N}%
 ( 0,\Omega_{m}^{-1} ) $, $0\leq m\leq m_{\ast}$ and $f_{m}$
denoting the corresponding joint density, and we define
%
\begin{equation}
\bar{\mathbb{P}}=\frac{1}{m_{\ast}}\sum_{m=1}^{m_{\ast}}
\mathbb {P}_{\Omega
_{m}}. \label{Pbar}
\end{equation}
For two distributions $\mathbb{P}$ and $\mathbb{Q}$ with densities $p$
and $%
q $ with respect to any common dominating measure $\mu$, we denote the
total variation affinity by $\Vert\mathbb{P}\wedge\mathbb{Q}\Vert
=\int
p\wedge q\,d\mu$.\vadjust{\goodbreak} The following lemma is a version of Le Cam's method; cf.
\citet{LeCam73}, \citet{Yu97}.\vadjust{\goodbreak}

\begin{lemma}
\label{Le Cam's lemma}Let $X^{(i)}$ be i.i.d. $\mathcal
{N}(0,\Omega
^{-1}) $, $i=1,2,\ldots,n$, with $\Omega\in\mathcal{G}_{0}$. Let
$\hat{%
\Omega}= ( \hat{\omega}_{kl} ) _{p\times p}$ be an estimator
of $%
\Omega_{m}= ( \omega_{kl}^{ ( m ) } ) _{p\times
p}$, then
\[
\sup_{0\leq m\leq m_{\ast}}\mathbb{P}_{\Omega_{m}} \biggl\{ \bigl\llvert
\hat{\omega}%
_{ij}-\omega_{ij}^{ ( m ) }\bigr
\rrvert >\frac{\alpha
}{2} \biggr\} \geq\frac{1}{2}\llVert
\mathbb{P}_{\Omega_{0}}\wedge\bar {\mathbb {P}}%
\rrVert ,
\]
where $\alpha=\inf_{1\leq m\leq m_{\ast}}\llvert \omega
_{ij}^{ (
m ) }-\omega_{ij}^{ ( 0 ) }\rrvert $.
\end{lemma}

\begin{pf*}{Proof of Theorem~\ref{sparseMinimaxlwSupport}} We shall
divide the
proof into three steps. Without loss of generality, consider only the
cases $%
 ( i,j ) = ( 1,1 ) $ and $ ( i,j )
= (
1,2 ) $. For the general case\vadjust{\goodbreak} $\omega_{ii}$ or $\omega_{ij}$
with $%
i\neq j$, we could always permute the coordinates and rearrange them to the
special case $\omega_{11}$ or $\omega_{12}$.

\begin{figure}

\includegraphics{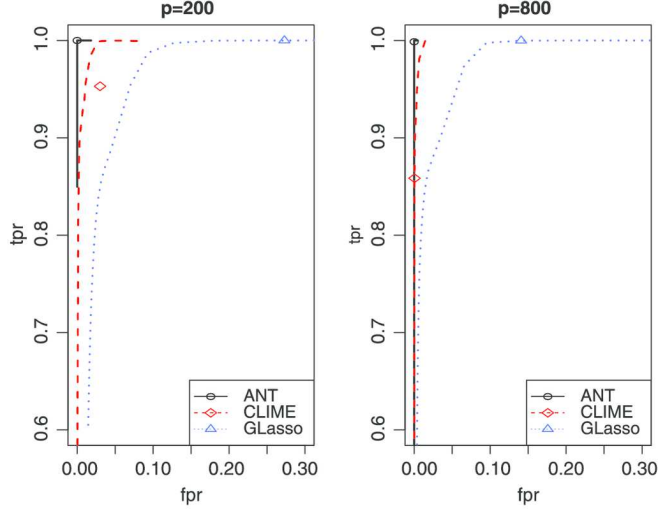}

\caption{The ROC curves.
Circle: ANT with the proposed thresholding.
Triangle: GLasso with penalty level by CV.
Diamond: CLIME with penalty level by CV.}
\label{fig-roc}
\end{figure}

\textit{Step} 1: \textit{Constructing the parameter set}. We first
define $%
\Omega_{0}$,
%
\begin{eqnarray}\label{SigmaZero}
\Sigma_{0}&=&\pmatrix{ 1 & b & 0 &
\cdots& 0
\vspace*{2pt}\cr
b & 1 & 0 & \cdots& 0
\vspace*{2pt}\cr
0 & 0 & 1 & \cdots& 0
\vspace*{2pt}\cr
\vdots& \vdots& \vdots& \ddots& \vdots
\vspace*{2pt}\cr
0 & 0 & 0 & 0 & 1}\quad\mbox{and}
\nonumber
\\[-8pt]
\\[-8pt]
\nonumber
\Omega_{0}&=&\Sigma_{0}^{-1}=\pmatrix{
\displaystyle\frac{1}{1-b^{2}} & \displaystyle\frac{-b}{1-b^{2}} & 0 &
\cdots& 0
\vspace*{2pt}\cr
\displaystyle\frac{-b}{1-b^{2}} & \displaystyle\frac{1}{1-b^{2}} & 0 & \cdots& 0
\vspace*{2pt}\cr
0 & 0 & 1 & \cdots& 0
\vspace*{2pt}\cr
\vdots& \vdots& \vdots& \ddots& \vdots
\vspace*{2pt}\cr
0 & 0 & 0 & 0 & 1};
\end{eqnarray}
that is, $\Sigma_{0}= ( \sigma_{kl}^{(0)} ) _{p\times p}$ is
a matrix
with all diagonal entries equal to 1, $\sigma_{12}^{(0)}=\sigma
_{21}^{(0)}=b$ and the rest all zeros. Here the constant $0<b<1$ is to be
determined later. For $\Omega_{m},1\leq m\leq m_{\ast}$, the construction
is as follows. Without loss of generality we assume $k_{n,p}\geq3$. Denote
by $\mathcal{H}$ the collection of all $p\times p$ symmetric matrices with
exactly $ ( k_{n,p}-2 ) $ elements equal to $1$ between the third
and the last elements on the first row (column) and the rest all zeros.
Define
%
\begin{equation}
\mathcal{G}_{0}= \bigl\{ \Omega\dvtx \Omega=\Omega_{0}
\mbox{ or }\Omega = ( \Sigma_{0}+aH ) ^{-1},\mbox{ for some }H
\in\mathcal{H} \bigr\} , \label{G0}
\end{equation}
where $a=\sqrt{\frac{\tau_{1}\log p}{n}}$ for some constant $\tau_{1}$
which is determined later. The cardinality of $\mathcal
{G}_{0}\setminus
 \{
\Omega_{0} \} $ is
\[
m^{\ast}=\mathrm{Card} ( \mathcal{G}_{0} ) -1=\mathrm {Card}
( \mathcal{H} ) =\pmatrix{p-2
\cr
k_{n,p}-2}.
\]
We pick the constant $b=\frac{1}{2} ( 1-1/M ) $ and
\[
0<\tau _{1}<\min \biggl\{ \frac{ ( 1-1/M )
^{2}-b^{2}}{C_{0}},\frac
{ (
1-b^{2} ) ^{2}}{2C_{0} ( 1+b^{2} ) },
\frac{ (
1-b^{2} ) ^{2}}{4\nu ( 1+b^{2} ) } \biggr\},
\]
and prove that $%
\mathcal{G}_{0}\subset\mathcal{G}_{0}(M,k_{n,p})$.

First we show that for all $\Omega_{i}$,
%
\begin{equation}
1/M\leq\lambda_{\min} ( \Omega_{i} ) <\lambda_{\max
}
( \Omega_{i} ) \leq M. \label{spectrum}
\end{equation}
For any matrix $\Omega_{m}$, $1\leq m\leq m_{\ast}$, some elementary
calculations yield that
\begin{eqnarray*}
\lambda_{1} \bigl( \Omega_{m}^{-1} \bigr) &=&1+
\sqrt{b^{2}+ ( k_{n,p}-2 ) a^{2}},\qquad
\lambda_{p} \bigl( \Omega_{m}^{-1} \bigr) =1-
\sqrt{%
b^{2}+ ( k_{n,p}-2 ) a^{2}},
\\
\lambda_{2} \bigl( \Omega_{m}^{-1} \bigr) &=&
\lambda_{3} \bigl( \Omega _{m}^{-1} \bigr) =
\cdots=\lambda_{p-1} \bigl( \Omega_{m}^{-1} \bigr)
=1.
\end{eqnarray*}
Since $b=\frac{1}{2} ( 1-1/M ) $ and $0<\tau_{1}<\frac
{ (
1-1/M ) ^{2}-b^{2}}{C_{0}}$, we have
%
\begin{eqnarray}\label{constant}
1-\sqrt{b^{2}+ ( k_{n,p}-2 ) a^{2}} &\geq&1-\sqrt
{b^{2}+\tau _{1}C_{0}}>1/M,
\nonumber
\\[-8pt]
\\[-8pt]
\nonumber
1+\sqrt{b^{2}+ ( k_{n,p}-2 ) a^{2}} &<&2-1/M<M,
\nonumber
\end{eqnarray}
which imply
\[
1/M\leq\lambda_{1}^{-1} \bigl( \Omega_{m}^{-1}
\bigr) =\lambda _{\min
} ( \Omega_{m} ) <
\lambda_{\max} ( \Omega_{m} ) =\lambda_{p}^{-1}
\bigl( \Omega_{m}^{-1} \bigr) \leq M.
\]
As for matrix $\Omega_{0}$, similarly we have
\begin{eqnarray*}
\lambda_{1} \bigl( \Omega_{0}^{-1} \bigr) &=&1+b,
\lambda_{p} \bigl( \Omega _{0}^{-1} \bigr) =1-b,
\\
\lambda_{2} \bigl( \Omega_{0}^{-1} \bigr) &=&
\lambda_{3} \bigl( \Omega _{0}^{-1} \bigr) =
\cdots=\lambda_{p-1} \bigl( \Omega_{0}^{-1} \bigr)
=1,
\end{eqnarray*}
and thus $1/M\leq\lambda_{\min} ( \Omega_{0} ) <\lambda
_{\max
} ( \Omega_{0} ) \leq M$ for the choice of $b=\frac
{1}{2} (
1-1/M ) $.

Now we show that the number of nonzero elements in $\Omega
_{m}$, $0\leq m\leq m_{\ast}$ is no more than $k_{n,p}$ per row/column.
From the construction of $\Omega_{m}^{-1}$, there exists some permutation
matrix $P_{\pi}$ such that $P_{\pi}\Omega_{m}^{-1}P_{\pi}^{T}$ is a
two-block diagonal matrix with dimensions $k_{n,p} $ and $%
 ( p-k_{n,p} ) $, of which the second block is an identity
matrix. Then $ ( P_{\pi}\Omega_{m}^{-1}P_{\pi}^{T} )
^{-1}=P_{\pi}\Omega_{m}P_{\pi}^{T}$ has the same blocking structure with
the first block of dimension $k_{n,p} $ and the second block
being an identity matrix. Thus the number of nonzero elements
is no more than $k_{n,p}$ per row/column for~$\Omega_{m}$. Therefore, we
have $\mathcal{G}_{0}\subset\mathcal{G}_{0}(M,k_{n,p})$ from equation
(\ref%
{spectrum}).

\textit{Step} 2: \textit{Bounding} $\bolds{\alpha}$. From the construction
of $\Omega_{m}^{-1}$ and the matrix inverse formula, we have that
for any precision matrix $\Omega_{m}$,
\[
\omega_{11}^{(m)}=\frac{1}{1-b^{2}- ( k_{n,p}-2 ) a^{2}}\quad \mbox{and}\quad
\omega_{12}^{(m)}=\frac{-b}{1-b^{2}- ( k_{n,p}-2 ) a^{2}}
\]
for $1\leq m\leq m_{\ast}$, and for the precision matrix $\Omega_{0}$,
\[
\omega_{11}^{(0)}=\frac{1}{1-b^{2}},\qquad\omega_{12}^{(0)}=
\frac{-b}{1-b^{2}}.
\]
Since $b^{2}+ ( k_{n,p}-2 ) a^{2}<$ $ ( 1-1/M )
^{2}<1$ in
equation (\ref{constant}), we have
%
\begin{eqnarray} \label{lossbd}
\inf_{1\leq m\leq m_{\ast}}\bigl\llvert \omega_{11}^{(m)}-
\omega _{11}^{(0)}\bigr\rrvert &=&\frac{ ( k_{n,p}-2 )
a^{2}}{ (
1-b^{2} )  ( 1-b^{2}- ( k_{n,p}-2 ) a^{2} )
}\geq
C_{3}k_{n,p}a^{2},
\nonumber
\\[-8pt]
\\[-8pt]
\nonumber
\inf_{1\leq m\leq m_{\ast}}\bigl\llvert \omega_{12}^{(m)}-
\omega _{12}^{(0)}\bigr\rrvert &=&\frac{b ( k_{n,p}-2 )
a^{2}}{ (
1-b^{2} )  ( 1-b^{2}- ( k_{n,p}-2 ) a^{2} )
}\geq
C_{4}k_{n,p}a^{2},
\nonumber
\end{eqnarray}
for some constants $C_{3},C_{4}>0$.

\textit{Step} 3: \textit{Bounding the affinity}. The following lemma is proved
in \citet{Ren13supp}.

\begin{lemma}
\label{affinity}Let $\bar{\mathbb{P}}$ be defined in (\ref{Pbar}).
We have
%
\begin{equation}
\llVert \mathbb{P}_{\Omega_{0}}\wedge\bar{\mathbb{P}}\rrVert \geq
C_{5} \label{affbd3}
\end{equation}
for some constant $C_{5}>0$.
\end{lemma}

Lemma~\ref{Le Cam's lemma}, together with equations (\ref{lossbd}),
(\ref%
{affbd3}) and $a=\sqrt{\frac{\tau_{1}\log p}{n}}$, imply
\begin{eqnarray*}
\sup_{0\leq m\leq m_{\ast}}\mathbb{P} \biggl\{ \bigl\llvert \hat {
\omega}%
_{11}-\omega_{11}^{ ( m ) }\bigr
\rrvert >\frac{1}{2}\cdot \frac
{%
C_{3}\tau_{1}k_{n,p}\log p}{n} \biggr\} &
\geq&C_{5}/2,
\\
\sup_{0\leq m\leq m_{\ast}}\mathbb{P} \biggl\{ \bigl\llvert \hat {
\omega}%
_{12}-\omega_{12}^{ ( m ) }\bigr
\rrvert >\frac{1}{2}\cdot \frac
{%
C_{4}\tau_{1}k_{n,p}\log p}{n} \biggr\} &
\geq&C_{5}/2,
\end{eqnarray*}
which match the lower bound in (\ref{sparserateSupport}) by setting $%
C_{1}=\min \{ C_{3}\tau_{1}/2,C_{4}\tau_{1}/2 \} $ and $%
c_{1}=C_{5}/2$.
\end{pf*}

\begin{remark}
Note that $|\!|\!|\Omega_{m}|\!|\!|_{1}$ is of the order $k_{n,p}\sqrt{\frac
{\log p}{%
n}}$, which implies $\frac{k_{n,p}\log p}{n}=k_{n,p}\sqrt{\frac{\log
p}{n}}%
\cdot\sqrt{\frac{\log p}{n}}\asymp|\!|\!|\Omega_{m}|\!|\!|_{1}\sqrt{\frac
{\log p}{%
n}}$. This observation partially explains why in the
literature we
need to assume the bounded matrix $\ell_{1}$ norm of
$\Omega
$ to
derive the lower bound rate $\sqrt{\frac{\log p}{n}}$. For the least
favorable parameter space, the matrix $\ell_{1}$ norm
of $\Omega$ cannot be avoided in the upper bound. However, the methodology proposed
in this paper improves the upper bounds in the literature by replacing the
matrix $\ell_{1}$ norm for every $\Omega$ by only
matrix $\ell_{1}$ norm bound of $\Omega$ in the least favorable
parameter space.
\end{remark}

%
\begin{supplement}[id=suppA]
\stitle{Supplement to ``Asymptotic normality and
optimalities in estimation of large Gaussian graphical
model''}
\slink[doi]{10.1214/14-AOS1286SUPP} 
\sdatatype{.pdf}
\sfilename{aos1286\_supp.pdf}
\sdescription{In this supplement we collect proofs of Theorems \ref
{Mainl0}--\ref{Main1} in Section~\ref{methodinf}, proofs of Theorems
\ref{support recovery}, \ref{Latent Graphical Model Inference} in
Section~\ref{app} and proofs of Theorems \ref{th-lasso}--\ref{Error
bound of LSE} as well as Proposition~\ref{prop-probab} in Section~\ref{regression-revisit}.}
\end{supplement}

\printaddresses

\begin{thebibliography}{47}

\bibitem[\protect\citeauthoryear{Antoniadis}{2010}]{Antoniadis10}
\begin{barticle}[mr]
\bauthor{\bsnm{Antoniadis},~\bfnm{Anestis}\binits{A.}}
(\byear{2010}).
\btitle{Comment: {$\ell_1$}-penalization for mixture regression models [MR2677722]}.
\bjournal{TEST}
\bvolume{19}
\bpages{257--258}.
\bid{doi={10.1007/s11749-010-0198-y}, issn={1133-0686}, mr={2677723}}
\end{barticle}
%

\bptok{imsref}%
\endbibitem

\bibitem[\protect\citeauthoryear{Belloni, Chernozhukov and Hansen}{2014}]{BelloniCH12}
\begin{barticle}[mr]
\bauthor{\bsnm{Belloni},~\bfnm{Alexandre}\binits{A.}},
\bauthor{\bsnm{Chernozhukov},~\bfnm{Victor}\binits{V.}} \AND
\bauthor{\bsnm{Hansen},~\bfnm{Christian}\binits{C.}}
(\byear{2014}).
\btitle{Inference on treatment effects after selection among high-dimensional controls}.
\bjournal{Rev. Econ. Stud.}
\bvolume{81}
\bpages{608--650}.
\bid{doi={10.1093/restud/rdt044}, issn={0034-6527}, mr={3207983}}
\end{barticle}
%

\bptok{imsref}%
\endbibitem

\bibitem[\protect\citeauthoryear{Belloni, Chernozhukov and Wang}{2011}]{BelloniCW11}
\begin{barticle}[mr]
\bauthor{\bsnm{Belloni},~\bfnm{A.}\binits{A.}},
\bauthor{\bsnm{Chernozhukov},~\bfnm{V.}\binits{V.}} \AND
\bauthor{\bsnm{Wang},~\bfnm{L.}\binits{L.}}
(\byear{2011}).
\btitle{Square-root lasso: Pivotal recovery of sparse signals via conic programming}.
\bjournal{Biometrika}
\bvolume{98}
\bpages{791--806}.
\bid{doi={10.1093/biomet/asr043}, issn={0006-3444}, mr={2860324}}
\end{barticle}
%

\bptok{imsref}%
\endbibitem

\bibitem[\protect\citeauthoryear{Bickel and Levina}{2008a}]{BL08A}
\begin{barticle}[mr]
\bauthor{\bsnm{Bickel},~\bfnm{Peter~J.}\binits{P.~J.}} \AND
\bauthor{\bsnm{Levina},~\bfnm{Elizaveta}\binits{E.}}
(\byear{2008}a).
\btitle{Regularized estimation of large covariance matrices}.
\bjournal{Ann. Statist.}
\bvolume{36}
\bpages{199--227}.
\bid{doi={10.1214/009053607000000758}, issn={0090-5364}, mr={2387969}}
\end{barticle}
%

\bptok{imsref}%
\endbibitem

\bibitem[\protect\citeauthoryear{Bickel and Levina}{2008b}]{BL08B}
\begin{barticle}[mr]
\bauthor{\bsnm{Bickel},~\bfnm{Peter~J.}\binits{P.~J.}} \AND
\bauthor{\bsnm{Levina},~\bfnm{Elizaveta}\binits{E.}}
(\byear{2008}b).
\btitle{Covariance regularization by thresholding}.
\bjournal{Ann. Statist.}
\bvolume{36}
\bpages{2577--2604}.
\bid{doi={10.1214/08-AOS600}, issn={0090-5364}, mr={2485008}}
\end{barticle}
%

\bptok{imsref}%
\endbibitem

\bibitem[\protect\citeauthoryear{Bickel, Ritov and Tsybakov}{2009}]{bickel2009simultaneous}
\begin{barticle}[mr]
\bauthor{\bsnm{Bickel},~\bfnm{Peter~J.}\binits{P.~J.}},
\bauthor{\bsnm{Ritov},~\bfnm{Ya'acov}\binits{Y.}} \AND
\bauthor{\bsnm{Tsybakov},~\bfnm{Alexandre~B.}\binits{A.~B.}}
(\byear{2009}).
\btitle{Simultaneous analysis of lasso and {D}antzig selector}.
\bjournal{Ann. Statist.}
\bvolume{37}
\bpages{1705--1732}.
\bid{doi={10.1214/08-AOS620}, issn={0090-5364}, mr={2533469}}
\end{barticle}
%

\bptok{imsref}%
\endbibitem

\bibitem[\protect\citeauthoryear{B{\"u}hlmann}{2013}]{Buhlmann12}
\begin{barticle}[mr]
\bauthor{\bsnm{B{\"u}hlmann},~\bfnm{Peter}\binits{P.}}
(\byear{2013}).
\btitle{Statistical significance in high-dimensional linear models}.
\bjournal{Bernoulli}
\bvolume{19}
\bpages{1212--1242}.
\bid{doi={10.3150/12-BEJSP11}, issn={1350-7265}, mr={3102549}}
\end{barticle}
%

\bptok{imsref}%
\endbibitem

\bibitem[\protect\citeauthoryear{Cai, Liu and Luo}{2011}]{CLL11}
\begin{barticle}[mr]
\bauthor{\bsnm{Cai},~\bfnm{Tony}\binits{T.}},
\bauthor{\bsnm{Liu},~\bfnm{Weidong}\binits{W.}} \AND
\bauthor{\bsnm{Luo},~\bfnm{Xi}\binits{X.}}
(\byear{2011}).
\btitle{A constrained {$\ell_1$} minimization approach to sparse precision matrix estimation}.
\bjournal{J. Amer. Statist. Assoc.}
\bvolume{106}
\bpages{594--607}.\vadjust{\goodbreak}
\bid{doi={10.1198/jasa.2011.tm10155}, issn={0162-1459}, mr={2847973}}
\end{barticle}
%

\bptok{imsref}%
\endbibitem

\bibitem[\protect\citeauthoryear{Cai, Liu and Zhou}{2012}]{CLZ12}
\begin{bmisc}[author]
\bauthor{\bsnm{Cai},~\bfnm{T~Tony}\binits{T.~T.}},
\bauthor{\bsnm{Liu},~\bfnm{Weidong}\binits{W.}} \AND
\bauthor{\bsnm{Zhou},~\bfnm{Harrison~H.}\binits{H.~H.}}
(\byear{2012}).
\bhowpublished{Estimating sparse precision matrix: {O}ptimal rates of convergence and adaptive estimation.
Preprint. Available at \arxivurl{arXiv:1212.2882}.}
\end{bmisc}
%

\bptok{imsref}%
\endbibitem

\bibitem[\protect\citeauthoryear{Cai, Zhang and Zhou}{2010}]{CZZH10}
\begin{barticle}[mr]
\bauthor{\bsnm{Cai},~\bfnm{T.~Tony}\binits{T.~T.}},
\bauthor{\bsnm{Zhang},~\bfnm{Cun-Hui}\binits{C.-H.}} \AND
\bauthor{\bsnm{Zhou},~\bfnm{Harrison~H.}\binits{H.~H.}}
(\byear{2010}).
\btitle{Optimal rates of convergence for covariance matrix estimation}.
\bjournal{Ann. Statist.}
\bvolume{38}
\bpages{2118--2144}.
\bid{doi={10.1214/09-AOS752}, issn={0090-5364}, mr={2676885}}
\end{barticle}
%

\bptok{imsref}%
\endbibitem

\bibitem[\protect\citeauthoryear{Cai and Zhou}{2012}]{CZh12Sparse}
\begin{barticle}[mr]
\bauthor{\bsnm{Cai},~\bfnm{T.~Tony}\binits{T.~T.}} \AND
\bauthor{\bsnm{Zhou},~\bfnm{Harrison~H.}\binits{H.~H.}}
(\byear{2012}).
\btitle{Optimal rates of convergence for sparse covariance matrix estimation}.
\bjournal{Ann. Statist.}
\bvolume{40}
\bpages{2389--2420}.
\bid{doi={10.1214/12-AOS998}, issn={0090-5364}, mr={3097607}}
\end{barticle}
%

\bptok{imsref}%
\endbibitem

\bibitem[\protect\citeauthoryear{Cand{\`e}s and Recht}{2009}]{CR09}
\begin{barticle}[mr]
\bauthor{\bsnm{Cand{\`e}s},~\bfnm{Emmanuel~J.}\binits{E.~J.}} \AND
\bauthor{\bsnm{Recht},~\bfnm{Benjamin}\binits{B.}}
(\byear{2009}).
\btitle{Exact matrix completion via convex optimization}.
\bjournal{Found. Comput. Math.}
\bvolume{9}
\bpages{717--772}.
\bid{doi={10.1007/s10208-009-9045-5}, issn={1615-3375}, mr={2565240}}
\end{barticle}
%

\bptok{imsref}%
\endbibitem

\bibitem[\protect\citeauthoryear{Chandrasekaran, Parrilo and Willsky}{2012}]{CPW12}
\begin{barticle}[mr]
\bauthor{\bsnm{Chandrasekaran},~\bfnm{Venkat}\binits{V.}},
\bauthor{\bsnm{Parrilo},~\bfnm{Pablo~A.}\binits{P.~A.}} \AND
\bauthor{\bsnm{Willsky},~\bfnm{Alan~S.}\binits{A.~S.}}
(\byear{2012}).
\btitle{Latent variable graphical model selection via convex optimization}.
\bjournal{Ann. Statist.}
\bvolume{40}
\bpages{1935--1967}.
\bid{doi={10.1214/11-AOS949}, issn={0090-5364}, mr={3059067}}
\end{barticle}
%

\bptok{imsref}%
\endbibitem

\bibitem[\protect\citeauthoryear{d'Aspremont, Banerjee and El~Ghaoui}{2008}]{d'ABE08}
\begin{barticle}[mr]
\bauthor{\bsnm{d'Aspremont},~\bfnm{Alexandre}\binits{A.}},
\bauthor{\bsnm{Banerjee},~\bfnm{Onureena}\binits{O.}} \AND
\bauthor{\bsnm{El Ghaoui},~\bfnm{Laurent}\binits{L.}}
(\byear{2008}).
\btitle{First-order methods for sparse covariance selection}.
\bjournal{SIAM J. Matrix Anal. Appl.}
\bvolume{30}
\bpages{56--66}.
\bid{doi={10.1137/060670985}, issn={0895-4798}, mr={2399568}}
\end{barticle}
%

\bptok{imsref}%
\endbibitem

\bibitem[\protect\citeauthoryear{El~Karoui}{2008}]{ELK08}
\begin{barticle}[mr]
\bauthor{\bsnm{El Karoui},~\bfnm{Noureddine}\binits{N.}}
(\byear{2008}).
\btitle{Operator norm consistent estimation of large-dimensional sparse covariance matrices}.
\bjournal{Ann. Statist.}
\bvolume{36}
\bpages{2717--2756}.
\bid{doi={10.1214/07-AOS559}, issn={0090-5364}, mr={2485011}}
\end{barticle}
%

\bptok{imsref}%
\endbibitem

\bibitem[\protect\citeauthoryear{Friedman, Hastie and Tibshirani}{2008}]{FHT08}
\begin{barticle}[pbm]
\bauthor{\bsnm{Friedman},~\bfnm{Jerome}\binits{J.}},
\bauthor{\bsnm{Hastie},~\bfnm{Trevor}\binits{T.}} \AND
\bauthor{\bsnm{Tibshirani},~\bfnm{Robert}\binits{R.}}
(\byear{2008}).
\btitle{Sparse inverse covariance estimation with the graphical lasso}.
\bjournal{Biostatistics}
\bvolume{9}
\bpages{432--441}.
\bid{doi={10.1093/biostatistics/kxm045}, issn={1468-4357}, mid={NIHMS248717}, pii={kxm045}, pmcid={3019769}, pmid={18079126}}
\end{barticle}
%

\bptok{imsref}%
\endbibitem

\bibitem[\protect\citeauthoryear{Horn and Johnson}{1990}]{MatAna}
\begin{bbook}[mr]
\bauthor{\bsnm{Horn},~\bfnm{Roger~A.}\binits{R.~A.}} \AND
\bauthor{\bsnm{Johnson},~\bfnm{Charles~R.}\binits{C.~R.}}
(\byear{1990}).
\btitle{Matrix Analysis}.
\bpublisher{Cambridge Univ. Press},
\blocation{Cambridge}.
\bid{mr={1084815}}
\end{bbook}
%

\bptok{imsref}%
\endbibitem

\bibitem[\protect\citeauthoryear{Javanmard and Montanari}{2014}]{JavanmardM13}
\begin{barticle}[mr]
\bauthor{\bsnm{Javanmard},~\bfnm{Adel}\binits{A.}} \AND
\bauthor{\bsnm{Montanari},~\bfnm{Andrea}\binits{A.}}
(\byear{2014}).
\btitle{Hypothesis testing in high-dimensional regression under the {G}aussian random design model: Asymptotic theory}.
\bjournal{IEEE Trans. Inform. Theory}
\bvolume{60}
\bpages{6522--6554}.
\bid{doi={10.1109/TIT.2014.2343629}, issn={0018-9448}, mr={3265038}}
\end{barticle}
%

\bptok{imsref}%
\endbibitem

\bibitem[\protect\citeauthoryear{Koltchinskii}{2009}]{Koltchinskii09}
\begin{barticle}[mr]
\bauthor{\bsnm{Koltchinskii},~\bfnm{Vladimir}\binits{V.}}
(\byear{2009}).
\btitle{The {D}antzig selector and sparsity oracle inequalities}.
\bjournal{Bernoulli}
\bvolume{15}
\bpages{799--828}.
\bid{doi={10.3150/09-BEJ187}, issn={1350-7265}, mr={2555200}}
\end{barticle}
%

\bptok{imsref}%
\endbibitem

\bibitem[\protect\citeauthoryear{Lam and Fan}{2009}]{LFAN09}
\begin{barticle}[mr]
\bauthor{\bsnm{Lam},~\bfnm{Clifford}\binits{C.}} \AND
\bauthor{\bsnm{Fan},~\bfnm{Jianqing}\binits{J.}}
(\byear{2009}).
\btitle{Sparsistency and rates of convergence in large covariance matrix estimation}.
\bjournal{Ann. Statist.}
\bvolume{37}
\bpages{4254--4278}.
\bid{doi={10.1214/09-AOS720}, issn={0090-5364}, mr={2572459}}
\end{barticle}
%

\bptok{imsref}%
\endbibitem

\bibitem[\protect\citeauthoryear{Lauritzen}{1996}]{GraMod}
\begin{bbook}[mr]
\bauthor{\bsnm{Lauritzen},~\bfnm{Steffen~L.}\binits{S.~L.}}
(\byear{1996}).
\btitle{Graphical Models}.
\bseries{Oxford Statistical Science Series}
\bvolume{17}.
\bpublisher{Oxford Univ. Press},
\blocation{New York}.
\bid{mr={1419991}}
\end{bbook}
%

\bptok{imsref}%
\endbibitem

\bibitem[\protect\citeauthoryear{Le Cam}{1973}]{LeCam73}
\begin{barticle}[mr]
\bauthor{\bsnm{Le Cam},~\bfnm{L.}\binits{L.}}
(\byear{1973}).
\btitle{Convergence of estimates under dimensionality restrictions}.
\bjournal{Ann. Statist.}
\bvolume{1}
\bpages{38--53}.
\bid{issn={0090-5364}, mr={0334381}}
\end{barticle}
%

\bptok{imsref}%
\endbibitem

\bibitem[\protect\citeauthoryear{Liu}{2013}]{Liu13}
\begin{barticle}[mr]
\bauthor{\bsnm{Liu},~\bfnm{Weidong}\binits{W.}}
(\byear{2013}).
\btitle{Gaussian graphical model estimation with false discovery rate control}.
\bjournal{Ann. Statist.}
\bvolume{41}
\bpages{2948--2978}.
\bid{doi={10.1214/13-AOS1169}, issn={0090-5364}, mr={3161453}}
\end{barticle}
%

\bptok{imsref}%
\endbibitem

\bibitem[\protect\citeauthoryear{Meinshausen and B{\"u}hlmann}{2006}]{MB06}
\begin{barticle}[mr]
\bauthor{\bsnm{Meinshausen},~\bfnm{Nicolai}\binits{N.}} \AND
\bauthor{\bsnm{B{\"u}hlmann},~\bfnm{Peter}\binits{P.}}
(\byear{2006}).
\btitle{High-dimensional graphs and variable selection with the lasso}.
\bjournal{Ann. Statist.}
\bvolume{34}
\bpages{1436--1462}.
\bid{doi={10.1214/009053606000000281}, issn={0090-5364}, mr={2278363}}
\end{barticle}
%

\bptok{imsref}%
\endbibitem

\bibitem[\protect\citeauthoryear{Pang, Liu and Vanderbei}{2014}]{JMLR:v15:pang14a}
\begin{barticle}[author]
\bauthor{\bsnm{Pang},~\bfnm{Haotian}\binits{H.}},
\bauthor{\bsnm{Liu},~\bfnm{Han}\binits{H.}} \AND
\bauthor{\bsnm{Vanderbei},~\bfnm{Robert}\binits{R.}}
(\byear{2014}).
\btitle{The FASTCLIME package for linear programming and large-scale precision matrix estimation in R}.
\bjournal{J. Mach. Learn. Res.}
\bvolume{15}
\bpages{489--493}.
\end{barticle}
%

\bptok{imsref}%
\endbibitem

\bibitem[\protect\citeauthoryear{Raskutti, Wainwright and Yu}{2010}]{raskutti2010restricted}
\begin{barticle}[mr]
\bauthor{\bsnm{Raskutti},~\bfnm{Garvesh}\binits{G.}},
\bauthor{\bsnm{Wainwright},~\bfnm{Martin~J.}\binits{M.~J.}} \AND
\bauthor{\bsnm{Yu},~\bfnm{Bin}\binits{B.}}
(\byear{2010}).
\btitle{Restricted eigenvalue properties for correlated {G}aussian designs}.
\bjournal{J. Mach. Learn. Res.}
\bvolume{11}
\bpages{2241--2259}.
\bid{issn={1532-4435}, mr={2719855}}
\end{barticle}
%

\bptok{imsref}%
\endbibitem

\bibitem[\protect\citeauthoryear{Ravikumar et~al.}{2011}]{RWRY08}
\begin{barticle}[mr]
\bauthor{\bsnm{Ravikumar},~\bfnm{Pradeep}\binits{P.}},
\bauthor{\bsnm{Wainwright},~\bfnm{Martin~J.}\binits{M.~J.}},
\bauthor{\bsnm{Raskutti},~\bfnm{Garvesh}\binits{G.}} \AND
\bauthor{\bsnm{Yu},~\bfnm{Bin}\binits{B.}}
(\byear{2011}).
\btitle{High-dimensional covariance estimation by minimizing {$\ell_1$}-penalized log-determinant divergence}.
\bjournal{Electron. J. Stat.}
\bvolume{5}
\bpages{935--980}.
\bid{doi={10.1214/11-EJS631}, issn={1935-7524}, mr={2836766}}
\end{barticle}
%

\bptok{imsref}%
\endbibitem

\bibitem[\protect\citeauthoryear{Ren and Zhou}{2012}]{RZ12}
\begin{barticle}[mr]
\bauthor{\bsnm{Ren},~\bfnm{Zhao}\binits{Z.}} \AND
\bauthor{\bsnm{Zhou},~\bfnm{Harrison~H.}\binits{H.~H.}}
(\byear{2012}).
\btitle{Discussion: {L}atent variable graphical model selection via convex optimization [MR3059067]}.
\bjournal{Ann. Statist.}
\bvolume{40}
\bpages{1989--1996}.
\bid{issn={0090-5364}, mr={3059072}}
\end{barticle}
%

\bptok{imsref}%
\endbibitem

\bibitem[\protect\citeauthoryear{Ren et al.}{2015}]{Ren13supp}
\begin{bmisc}[author]
\bauthor{\bsnm{Ren},~\bfnm{Zhao}\binits{Z.}},
\bauthor{\bsnm{Sun},~\bfnm{Tingni}\binits{T.}},
\bauthor{\bsnm{Zhang},~\bfnm{Cun-Hui}\binits{C.-H.}} \AND
\bauthor{\bsnm{Zhou},~\bfnm{Harrison~H.}\binits{H.~H.}}
(\byear{2015}).
\bhowpublished{Supplement to ``Asymptotic normality and optimalities in estimation of large
Gaussian graphical models.''
DOI:\doiurl{10.1214/14-AOS1286SUPP}}.
\bptok{imsref}%
\end{bmisc}
\endbibitem

\bibitem[\protect\citeauthoryear{Rothman et~al.}{2008}]{RBLZ08}
\begin{barticle}[mr]
\bauthor{\bsnm{Rothman},~\bfnm{Adam~J.}\binits{A.~J.}},
\bauthor{\bsnm{Bickel},~\bfnm{Peter~J.}\binits{P.~J.}},
\bauthor{\bsnm{Levina},~\bfnm{Elizaveta}\binits{E.}} \AND
\bauthor{\bsnm{Zhu},~\bfnm{Ji}\binits{J.}}
(\byear{2008}).
\btitle{Sparse permutation invariant covariance estimation}.
\bjournal{Electron. J. Stat.}
\bvolume{2}
\bpages{494--515}.
\bid{doi={10.1214/08-EJS176}, issn={1935-7524}, mr={2417391}}
\end{barticle}
%

\bptok{imsref}%
\endbibitem

\bibitem[\protect\citeauthoryear{St{\"a}dler, B{\"u}hlmann and van~de Geer}{2010}]{StadlerBG10}
\begin{barticle}[mr]
\bauthor{\bsnm{St{\"a}dler},~\bfnm{Nicolas}\binits{N.}},
\bauthor{\bsnm{B{\"u}hlmann},~\bfnm{Peter}\binits{P.}} \AND
\bauthor{\bsnm{van~de Geer},~\bfnm{Sara}\binits{S.}}
(\byear{2010}).
\btitle{{$\ell_1$}-penalization for mixture regression models}.
\bjournal{TEST}
\bvolume{19}
\bpages{209--256}.
\bid{doi={10.1007/s11749-010-0197-z}, issn={1133-0686}, mr={2677722}}
\bptnote{check related}%
\end{barticle}
%

\bptok{imsref}%
\endbibitem

\bibitem[\protect\citeauthoryear{Sun and Zhang}{2010}]{SunZ10}
\begin{barticle}[mr]
\bauthor{\bsnm{Sun},~\bfnm{Tingni}\binits{T.}} \AND
\bauthor{\bsnm{Zhang},~\bfnm{Cun-Hui}\binits{C.-H.}}
(\byear{2010}).
\btitle{Comment: {$\ell_1$}-penalization for mixture regression models [MR2677722]}.
\bjournal{TEST}
\bvolume{19}
\bpages{270--275}.
\bid{doi={10.1007/s11749-010-0201-7}, issn={1133-0686}, mr={2677726}}
\end{barticle}
%

\bptok{imsref}%
\endbibitem

\bibitem[\protect\citeauthoryear{Sun and Zhang}{2012a}]{SZH12}
\begin{barticle}[mr]
\bauthor{\bsnm{Sun},~\bfnm{Tingni}\binits{T.}} \AND
\bauthor{\bsnm{Zhang},~\bfnm{Cun-Hui}\binits{C.-H.}}
(\byear{2012}a).
\btitle{Scaled sparse linear regression}.
\bjournal{Biometrika}
\bvolume{99}
\bpages{879--898}.
\bid{doi={10.1093/biomet/ass043}, issn={0006-3444}, mr={2999166}}
\end{barticle}
%

\bptok{imsref}%
\endbibitem

\bibitem[\protect\citeauthoryear{Sun and Zhang}{2012b}]{SunZhang2012SinicaComment}
\begin{barticle}[mr]
\bauthor{\bsnm{Sun},~\bfnm{Tingni}\binits{T.}} \AND
\bauthor{\bsnm{Zhang},~\bfnm{Cun-Hui}\binits{C.-H.}}
(\byear{2012}b).
\btitle{Comment: ``{M}inimax estimation of large covariance matrices under {$\ell_1$}-norm'' [MR3027084]}.
\bjournal{Statist. Sinica}
\bvolume{22}
\bpages{1354--1358}.
\bid{issn={1017-0405}, mr={3027086}}
\end{barticle}
%

\bptok{imsref}%
\endbibitem

\bibitem[\protect\citeauthoryear{Sun and Zhang}{2013}]{SZH122}
\begin{barticle}[mr]
\bauthor{\bsnm{Sun},~\bfnm{Tingni}\binits{T.}} \AND
\bauthor{\bsnm{Zhang},~\bfnm{Cun-Hui}\binits{C.-H.}}
(\byear{2013}).
\btitle{Sparse matrix inversion with scaled lasso}.
\bjournal{J. Mach. Learn. Res.}
\bvolume{14}
\bpages{3385--3418}.
\bid{issn={1532-4435}, mr={3144466}}
\end{barticle}
%

\bptok{imsref}%
\endbibitem

\bibitem[\protect\citeauthoryear{Thorin}{1948}]{thorin1948convexity}
\begin{barticle}[mr]
\bauthor{\bsnm{Thorin},~\bfnm{G.~O.}\binits{G.~O.}}
(\byear{1948}).
\btitle{Convexity theorems generalizing those of M. {R}iesz and {H}adamard with some applications}.
\bjournal{Comm. Sem. Math. Univ. Lund [Medd. Lunds Univ. Mat. Sem.]}
\bvolume{9}
\bpages{1--58}.
\bid{mr={0025529}}
\end{barticle}
%

\bptok{imsref}%
\endbibitem

\bibitem[\protect\citeauthoryear{van~de Geer and B{\"u}hlmann}{2009}]{van2009conditions}
\begin{barticle}[mr]
\bauthor{\bsnm{van~de Geer},~\bfnm{Sara~A.}\binits{S.~A.}} \AND
\bauthor{\bsnm{B{\"u}hlmann},~\bfnm{Peter}\binits{P.}}
(\byear{2009}).
\btitle{On the conditions used to prove oracle results for the {L}asso}.
\bjournal{Electron. J. Stat.}
\bvolume{3}
\bpages{1360--1392}.
\bid{doi={10.1214/09-EJS506}, issn={1935-7524}, mr={2576316}}
\end{barticle}
%

\bptok{imsref}%
\endbibitem

\bibitem[\protect\citeauthoryear{van~de Geer et~al.}{2014}]{GeerBR13}
\begin{barticle}[mr]
\bauthor{\bsnm{van~de Geer},~\bfnm{Sara}\binits{S.}},
\bauthor{\bsnm{B{\"u}hlmann},~\bfnm{Peter}\binits{P.}},
\bauthor{\bsnm{Ritov},~\bfnm{Ya'acov}\binits{Y.}} \AND
\bauthor{\bsnm{Dezeure},~\bfnm{Ruben}\binits{R.}}
(\byear{2014}).
\btitle{On asymptotically optimal confidence regions and tests for high-dimensional models}.
\bjournal{Ann. Statist.}
\bvolume{42}
\bpages{1166--1202}.
\bid{doi={10.1214/14-AOS1221}, issn={0090-5364}, mr={3224285}}
\end{barticle}
%

\bptok{imsref}%
\endbibitem

\bibitem[\protect\citeauthoryear{Ye and Zhang}{2010}]{YZH10}
\begin{barticle}[mr]
\bauthor{\bsnm{Ye},~\bfnm{Fei}\binits{F.}} \AND
\bauthor{\bsnm{Zhang},~\bfnm{Cun-Hui}\binits{C.-H.}}
(\byear{2010}).
\btitle{Rate minimaxity of the {L}asso and {D}antzig selector for the {$\ell_q$} loss in {$\ell_r$} balls}.
\bjournal{J. Mach. Learn. Res.}
\bvolume{11}
\bpages{3519--3540}.
\bid{issn={1532-4435}, mr={2756192}}
\end{barticle}
%

\bptok{imsref}%
\endbibitem

\bibitem[\protect\citeauthoryear{Yu}{1997}]{Yu97}
\begin{bmisc}[author]
\bauthor{\bsnm{Yu},~\bfnm{B.}\binits{B.}}
(\byear{1997}).
\bhowpublished{Assouad, {F}ano, and {L}e {C}am.
In \textit{Festschrift for Lucien Le Cam}
423--435. Springer, New York}.
\end{bmisc}
%

\bptok{imsref}%
\endbibitem

\bibitem[\protect\citeauthoryear{Yuan}{2010}]{Yuan10}
\begin{barticle}[mr]
\bauthor{\bsnm{Yuan},~\bfnm{Ming}\binits{M.}}
(\byear{2010}).
\btitle{High dimensional inverse covariance matrix estimation via linear programming}.
\bjournal{J. Mach. Learn. Res.}
\bvolume{11}
\bpages{2261--2286}.
\bid{issn={1532-4435}, mr={2719856}}
\end{barticle}
%

\bptok{imsref}%
\endbibitem

\bibitem[\protect\citeauthoryear{Yuan and Lin}{2007}]{YuanLin07}
\begin{barticle}[mr]
\bauthor{\bsnm{Yuan},~\bfnm{Ming}\binits{M.}} \AND
\bauthor{\bsnm{Lin},~\bfnm{Yi}\binits{Y.}}
(\byear{2007}).
\btitle{Model selection and estimation in the {G}aussian graphical model}.
\bjournal{Biometrika}
\bvolume{94}
\bpages{19--35}.
\bid{doi={10.1093/biomet/asm018}, issn={0006-3444}, mr={2367824}}
\end{barticle}
%

\bptok{imsref}%
\endbibitem

\bibitem[\protect\citeauthoryear{Zhang}{2009}]{Zhang09-l1}
\begin{barticle}[mr]
\bauthor{\bsnm{Zhang},~\bfnm{Tong}\binits{T.}}
(\byear{2009}).
\btitle{Some sharp performance bounds for least squares regression with {$L_1$} regularization}.
\bjournal{Ann. Statist.}
\bvolume{37}
\bpages{2109--2144}.
\bid{doi={10.1214/08-AOS659}, issn={0090-5364}, mr={2543687}}
\end{barticle}
%

\bptok{imsref}%
\endbibitem

\bibitem[\protect\citeauthoryear{Zhang}{2011}]{ZhangOber11}
\begin{bmisc}[author]
\bauthor{\bsnm{Zhang},~\bfnm{Cun-Hui}\binits{C.-H.}}
(\byear{2011}).
\bhowpublished{Statistical inference for high-dimensional data.
In \textit{Mathematisches Forschungsinstitut Oberwolfach:
Very High Dimensional Semiparametric Models}. Report No.~48/2011
28--31}.
\end{bmisc}
%

\bptok{imsref}%
\endbibitem

\bibitem[\protect\citeauthoryear{Zhang and Huang}{2008}]{ZhangH08}
\begin{barticle}[mr]
\bauthor{\bsnm{Zhang},~\bfnm{Cun-Hui}\binits{C.-H.}} \AND
\bauthor{\bsnm{Huang},~\bfnm{Jian}\binits{J.}}
(\byear{2008}).
\btitle{The sparsity and bias of the LASSO selection in high-dimensional linear regression}.
\bjournal{Ann. Statist.}
\bvolume{36}
\bpages{1567--1594}.
\bid{doi={10.1214/07-AOS520}, issn={0090-5364}, mr={2435448}}
\end{barticle}
%

\bptok{imsref}%
\endbibitem

\bibitem[\protect\citeauthoryear{Zhang and Zhang}{2012}]{ZhangZ12StatSc}
\begin{barticle}[mr]
\bauthor{\bsnm{Zhang},~\bfnm{Cun-Hui}\binits{C.-H.}} \AND
\bauthor{\bsnm{Zhang},~\bfnm{Tong}\binits{T.}}
(\byear{2012}).
\btitle{A general theory of concave regularization for high-dimensional sparse estimation problems}.
\bjournal{Statist. Sci.}
\bvolume{27}
\bpages{576--593}.
\bid{doi={10.1214/12-STS399}, issn={0883-4237}, mr={3025135}}
\end{barticle}
%

\bptok{imsref}%
\endbibitem

\bibitem[\protect\citeauthoryear{Zhang and Zhang}{2014}]{zhang2014confidence}
\begin{barticle}[mr]
\bauthor{\bsnm{Zhang},~\bfnm{Cun-Hui}\binits{C.-H.}} \AND
\bauthor{\bsnm{Zhang},~\bfnm{Stephanie~S.}\binits{S.~S.}}
(\byear{2014}).
\btitle{Confidence intervals for low dimensional parameters in high dimensional linear models}.
\bjournal{J. R. Stat. Soc. Ser. B Stat. Methodol.}
\bvolume{76}
\bpages{217--242}.
\bid{doi={10.1111/rssb.12026}, issn={1369-7412}, mr={3153940}}
\end{barticle}
%

\bptok{imsref}%
\endbibitem
\end{thebibliography}
\end{document}